\documentclass[preprint,12pt]{elsarticle}

\usepackage{amsmath,amssymb,eepic}

\renewcommand{\le}{\leqslant}
\renewcommand{\ge}{\geqslant}
\renewcommand{\sec}{\cap}
\renewcommand{\phi}{\varphi}
\newcommand{\TO}{\rightrightarrows}

\newcommand{\union}{\cup}

\makeatletter\let\ellipsearc\@arc\makeatother
\newcommand{\id}{\mathrm{id}}
\newcommand{\0}{\emptyset}
\newcommand{\1}{\mathbf{1}}
\newcommand\UNION\bigcup

\newcommand{\MM}{\mathcal{M}}

\usepackage{theorem}

\theorembodyfont{\slshape}

  \newtheorem{THM}{Theorem}[section]
  \newtheorem{LEM}[THM]{Lemma}
  
  \newtheorem{COR}[THM]{Corollary}

\theorembodyfont{\rmfamily}

  \newtheorem{DEF}[THM]{Definition}

\newif\ifQEDsign
\newcommand{\QED}{\global\QEDsigntrue\hfill$\square$}

\newenvironment{PROOF}%
    {\par\noindent\textit{Proof.}\global\QEDsignfalse}%
    {\ifQEDsign\else\QED\fi\par\bigskip\par}

\journal{Discrete Mathematics}

\begin{document}

\begin{frontmatter}

\title{On finite reflexive homomorphism-homogeneous\\binary relational systems\tnoteref{THX}}
\tnotetext[THX]{Supported by the Grant No 144017 of the Ministry of Science of the republic of Serbia}

\author{Dragan Ma\v sulovi\'c, Rajko Nenadov, Nemanja \v Skori\'c}
\ead{$\{$dragan.masulovic, rajko.nenadov, nemanja.skoric$\}$@dmi.uns.ac.rs}
\address{Department of Mathematics and Informatics, University of Novi Sad\\
        Trg Dositeja Obradovi\'ca 4, 21000 Novi Sad, Serbia}

\begin{abstract}
  A structure is called homogeneous if every isomorphism between finitely induced
  substructures of the structure extends to an automorphism of the structure.
  Recently, P.~J.~Cameron and J.~Ne\v set\v ril introduced a relaxed version
  of homogeneity: we say that a structure is homo\-mor\-phism-homogeneous if every
  homomorphism between finitely induced substructures of the structure extends to an endomorphism
  of the structure.

  In this paper we consider finite homomorphism-homogeneous relational systems with one reflexive binary relation.
  We show that for a large part of such relational systems (bidirectionally connected digraphs;
  a digraph is bidirectionally connected if each of its connected components can be traversed
  by $\rightleftarrows$-paths) the problem of deciding whether the system is homomorphism-homogeneous is
  coNP-complete. Consequently, for this class of relational systems we cannot hope for a description
  involving a catalogue, where by a catalogue we understand a finite list of polynomially decidable classes of
  structures. On the other hand, in case of bidirectionally disconnected digraphs we present
  the full characterization. Our main result states that if a digraph is bidirectionally
  disconnected, then it is homomorphism-homogeneous if and only if it is either a finite homomorphism-homogeneous
  quasiorder, or an inflation of a homomorphism-homogeneous digraph with involution
  (a peculiar class of digraphs introduced later in the paper), or an
  inflation of a digraph whose only connected components are $C_3^\circ$ and~$\1^\circ$.
\end{abstract}


\end{frontmatter}

\section{Introduction}

A structure is \emph{homogeneous} if every isomorphism between finitely induced
substructures of the structure extends to an automorphism of the structure.
In their recent paper \cite{cameron-nesetril} the authors discuss
a generalization of homogeneity to various types of morphisms between structures,
and in particular introduce the notion of homomorphism-homogeneous structures:

\begin{DEF}[Cameron, Ne\v set\v ril \cite{cameron-nesetril}]\label{nans.def.1}
  A structure is called \emph{homo\-mor\-phism-homogeneous} if every homomorphism between
  finitely induced substructures of the structure extends to an endomorphism of the structure.
\end{DEF}

This paper grew out of the authors' intention to characterize all
finite homomorphism-homogeneous relational systems with one reflexive binary relation
(binary relational systems). However, the complete
characterization of such relational systems turns out to be rather involved since the presence of
loops allows homomorphisms to spread their wings. What makes the problem unsolvable in general
is a result presented in \cite{rus-sch} where the authors show
that the problem of deciding whether a finite graph with loops allowed is homomorphism-homogeneous
is coNP-complete.

After the introductory Section~\ref{nans.sec.prelim},
in Section~\ref{nans.sec.theta-conn} we adapt the argument of \cite{rus-sch} to show that
the same holds even for bidirectionally connected improper digraphs
(a digraph is bidirectionally connected if each of its connected components can be traversed
by $\rightleftarrows$-paths; it is improper if it contains both edges of the form $\rightleftarrows$
and of the form~$\rightarrow$). The fact that deciding homomorphism-homogeneity
is computationally hard for bidirectionally connected digraphs means that for this class of digraphs
we cannot hope for a full description that involves a catalogue, where by a catalogue we
understand a finite list of polynomially decidable classes of structures.
We then turn to the classification of bidirectionally disconnected systems, which heavily relies
on a peculiar class of digraphs we refer to as digraphs with involution. Section~\ref{nans.sec.dinv}
is devoted to the classification of homomorphism-homogeneous digraphs in that class.
A rather long Section~\ref{nans.sec.bi-disconn} concludes the paper and classifies all finite reflexive
homomorphism-homogeneous bidirectionally disconnected systems.
Our main result is Corollary~\ref{nans.cor.MAIN} which states that if a digraph is bidirectionally
disconnected, then it is homomorphism-homogeneous if and only if it is either a finite homomorphism-homogeneous
quasiorder, or an inflation of a homomorphism-homogeneous digraph with involution, or an
inflation of a digraph whose only connected components are $C_3^\circ$ and~$\1^\circ$.

\section{Preliminaries}
\label{nans.sec.prelim}

A binary relational system is an ordered pair $(V, E)$ where $E \subseteq V^2$ is a binary
relation on $V$. A binary relational system $(V, E)$ is \emph{reflexive} if $(x, x) \in E$ for
all $x \in V$, \emph{irreflexive} if $(x, x) \notin E$ for all $x \in V$, \emph{symmetric}
if $(x, y) \in E$ implies $(y, x) \in E$ for all $x, y \in V$ and \emph{antisymmetric}
if $(x, y) \in E$ implies $(y, x) \notin E$ for all distinct $x, y \in V$.

Binary relational systems can be thought of in terms of digraphs (hence the notation $(V, E)$).
Then $V$ is the set of \emph{vertices} and $E$ is the set of \emph{edges} of the binary relational
system/digraph $(V, E)$. Edges of the form $(x, x)$ are called \emph{loops}. If $(x, x) \in E$
we also say that \emph{$x$ has a loop}. Instead of $(x, y) \in E$ we often
write $x \to y$ and say that $x$ \emph{dominates} $y$, or that $y$ \emph{is dominated by} $x$.
By $x \sim y$ we denote that $x \to y$ or $y \to x$, while
$x \rightleftarrows y$ denotes that $x \to y$ and $y \to x$. If $x \rightleftarrows y$,
we say that $x$ and $y$ form a \emph{double edge}. We shall also say that
a vertex $x$ is \emph{incident with a double edge} if there is a vertex $y \ne x$ such that
$x \rightleftarrows y$.

Digraphs $(V, E)$ where $E$ is a symmetric binary relation on $V$ are usually referred to as \emph{graphs}.
\emph{Proper digraphs} are digraphs $(V, E)$ where $E$ is an antisymmetric binary relation.
In this paper, digraphs $(V, E)$ where $E$ is neither antisymmetric nor symmetric will be referred to
as \emph{improper digraphs}. In an improper digraph there exists a pair of distinct
vertices $x$ and $y$ such that $x \rightleftarrows y$ and another pair of distinct vertices
$u$ and $v$ such that $u \to v$ and $v \not\to u$.

If $X, Y \subseteq V$ are nonempty
subsets of $V$ then $X \to Y$ means that $x \to y$ for some $x \in X$ and some $y \in Y$.
By $X \sim Y$ we denote that $X \to Y$ or $Y \to X$, while $X \rightleftarrows Y$ denotes that $X \to Y$ and
$Y \to X$. Moreover, $X \TO Y$ stands for $x \to y$ for all $x \in X$ and $y \in Y$.
Instead of $\{x\} \TO Y$ and $X \TO \{y\}$ we write $x \TO Y$ and $X \TO y$,
respectively, and analogously for $x \to Y$, $X \to y$ and $x \rightleftarrows Y$.
Let $r, s, t, u \in V(D)$ be vertices of a digraph $D$. We write $\{r, s\} \bowtie \{t, u\}$ to denote that
$r \to t \to s \to u \to r$ or $r \to u \to s \to t \to r$.

A digraph $D' = (V', E')$ is a \emph{subdigraph} of a digraph $D = (V, E)$ if
$V' \subseteq V$ and $E' \subseteq E$. We write $D' \le D$ to denote that
$D'$ is isomorphic to a subdigraph of~$D$.
For $\0 \ne W \subseteq V$ by $D[W]$ we denote the digraph $(W, E \sec W^2)$
which we refer to as the \emph{subdigraph of $D$ induced by $W$}.

Vertices $x$ and $y$ are \emph{connected in $D$} if there exists a sequence of vertices
$z_1, \ldots, z_k \in V$ such that $x = z_1 \sim \ldots \sim z_k = y$. A digraph $D$ is
\emph{weakly connected} if each pair of distinct vertices of $D$ is connected in~$D$.
A digraph $D$ is \emph{disconnected} if it is not weakly connected.
A \emph{connected component of $D$} is a maximal set $S \subseteq V$ such that
$D[S]$ is weakly connected. The number of connected components of $D$ will be denoted by $\omega(D)$.

Vertices $x$ and $y$ are \emph{bidirectionally connected in $D$} if there exists a sequence of vertices
$z_1, \ldots, z_k \in V$ such that $x = z_1 \rightleftarrows \ldots \rightleftarrows z_k = y$.
Define a binary relation $\theta(D)$ on $V(D)$ as follows:
$(x, y) \in \theta(D)$ if and only if $x = y$ or $x$ and $y$ are bidirectionally connected.
Clearly, $\theta(D)$ is an equivalence relation on $V(D)$ and $\omega(D) \le |V(D) / \theta(D)|$.
We say that a digraph $D$ is \emph{bidirectionally connected} if $\omega(D) = |V(D) / \theta(D)|$,
and that it is \emph{bidirectionally disconnected} if $\omega(D) < |V(D) / \theta(D)|$.
Note that a bidirectionally connected digraph need not be connected, and that a bidirectionally disconnected digraph
need not be disconnected; a digraph $D$ is bidirectionally connected if every connected component of $D$ contains
precisely one $\theta(D)$-class, while it is bidirectionally disconnected if there exists a connected component of
$D$ which consists of at least two $\theta(D)$-classes. In particular, every proper digraph with at least
two vertices is bidirectionally disconnected, and every graph
(even a disconnected one) is bidirectionally connected.

Let $K_n$ denote the complete irreflexive graph on $n$ vertices, and let
$K_n^\circ$ denote the complete reflexive graph on $n$ vertices.
Let $\1$ denote the trivial digraph with only one vertex and no edges,
and let $\1^\circ$ denote the digraph with only one vertex with a loop.
An \emph{oriented cycle with $n$ vertices} is a digraph $C_n$ whose vertices are
$1$, $2$, \ldots, $n$, $n \ge 3$, and whose only edges are $1 \to 2 \to \ldots \to n \to 1$.

For digraphs $D_1 = (V_1, E_1)$ and $D_2 = (V_2, E_2)$, by $D_1 + D_2$ we denote the
\emph{disjoint union} of $D_1$ and $D_2$. We assume that $D + O = O + D = D$, where $O = (\0, \0)$
denotes the \emph{empty digraph}.
The disjoint union $\underbrace{D + \ldots + D}_k$ consisting of $k \ge 1$ copies of $D$ will be abbreviated
to $k \cdot D$. Moreover, we let $0 \cdot D = O$.

Let $D_1 = (V_1, E_1)$ and $D_2 = (V_2, E_2)$ be digraphs. We say that
$f : V_1 \to V_2$ is a \emph{homomorphism} between $D_1$ and $D_2$ and write $f : D_1 \to D_2$
if
$$
  x \to y \text{ implies } f(x) \to f(y), \text{ for all } x, y \in V_1.
$$
An \emph{endomorphsim} is a homomorphism from $D$ into itself.
A mapping $f : V_1 \to V_2$ is an \emph{isomorphism} between $D_1$ and $D_2$ if
$f$ is bijective and
$$
  x \to y \text{ if and only if } f(x) \to f(y), \text{ for all } x, y \in V_1.
$$
Digraphs $D_1$ and $D_2$ are \emph{isomorphic} if there is an isomorphism between them.
We write $D_1 \cong D_2$. An \emph{automorphsim} is an isomorphism from $D$ onto itself.

A digraph $D$ is \emph{homomorphism-homogeneous} if every homomorphism $f : W_1 \to W_2$ between
finitely induced subdigraphs of $D$ extends to an endomorphism of~$D$
(see Definition~\ref{nans.def.1}).

For digraphs $D_1 = (V_1, E_1)$ and $D_2 = (V_2, E_2)$ we write $D_1 \Rightarrow D_2$ to denote that
every homomorphism $f : D_1[U] \to D_2[W]$, where $\0 \ne U \subseteq V_1$ and
$\0 \ne W \subseteq V_2$, extends to a
homomorphism $f^* : D_1 \to D_2$. Clearly, a digraph $D$ is homomorphism-homogeneous if and only if
$D \Rightarrow D$. The following two statements are obvious:

\begin{LEM}\label{hhdl.lem.disconn-iff}
  Let $D$ be a digraph. Then $D$ is homomorphism-homogeneous if and only if
  $D[S] \Rightarrow D[S']$ for every pair of (not necessarily distinct) connected components $S$, $S'$ of $D$.
\end{LEM}

\begin{LEM}\label{nans.lem.fS-S1}
  Let $D$ be an improper digraph and let $f$ be an endomorphism of $D$. Then
  \begin{itemize}
  \item
    for every connected component $S$ of $D$ there exists a connected component $S'$ of $D$
    such that $f(S) \subseteq S'$;
  \item
    for every $S \in V(D) / \theta(D)$ there exists an $S' \in V(D) / \theta(D)$ such that
    $f(S) \subseteq S'$.
  \end{itemize}
\end{LEM}

A digraph $D = (V, E)$ is \emph{transitive} if $x \to y \to z$ implies $x \to z$
for all $x, y, z \in V$.
Transitive reflexive proper digraphs are usually referred to as \emph{partially ordered sets}.
Recall that a mapping $f : A \to B$ is a \emph{homomorphism} between partially ordered sets $(A, \le)$
and $(B, \le)$ if
$$
  x \le y \text{ implies } f(x) \le f(y), \text{ for all } x, y \in A.
$$
It is clear that a mapping is a homomorphism between two partially ordered sets in the above sense
if and only if the mapping is a homomorphism between the corresponding digraphs.
Therefore, a paritally ordered set is homomorphism-homogeneous as a partially ordered set if and
only if it is homomorphism-homogeneous as a digraph.

\begin{THM}[\cite{masulovic-po}]
  \label{nans.thm.hh-poset}
  A finite partially ordered set $(A, \le)$ is homomorphism-homo\-gen\-eous
  if and only if
  \begin{itemize}
  \item[(1)]
    every connected component of $(A, \le)$ is a chain (this case includes anti-chains);
  \item[(2)]
    $(A, \le)$ is a tree, where a tree is a connected partially ordered set whose every up-set
    $\uparrow x$ is a chain;
  \item[(3)]
    $(A, \le)$ is a dual tree, where a dual tree is a connected partially ordered set whose every down-set
    $\downarrow x$ is a chain;
  \item[(4)]
    $(A, \le)$ splits into a tree and a dual tree in the following sense:
    there exists a partition $\{I, F\}$ of $A$ such that
    \begin{itemize}
    \item[$(i)$]
      $I$ is an ideal in $A$ and a tree,
    \item[$(ii)$]
      $F$ is a filter in $A$ and a dual tree, and
    \item[$(iii)$]
      $\forall x \in I \; \exists y \in F \; (x \le y)$ and
      $\forall y \in F \; \exists x \in I \; (x \le y)$;
    \end{itemize}
  \item[(5)]
    $(A, \le)$ is a lattice.
  \end{itemize}
\end{THM}

Reflexive finite homomorphism-homogeneous proper digraphs were characterized in
Theorem 3.10 of \cite{masul-hhd}:

\begin{THM}[\cite{masul-hhd}]
  \label{nans.thm.hh-digraph}
  Let $D$ be a reflexive proper digraph.
  Then $D$ is homo\-morph\-ism-homogeneous if and only if $D$ is one of the following digraphs:
  \begin{enumerate}
  \item
    $k \cdot C_3^\circ + l \cdot \1^\circ$ for some $k, l \ge 0$ such that
    $k + l \ge 1$;
  \item
    a finite homomorphism-homogeneous partially ordered set (see
    Theorem~\ref{nans.thm.hh-poset}).
  \end{enumerate}
\end{THM}

\section{Bidirectionally connected systems}
\label{nans.sec.theta-conn}

Rusinov and Schweitzer have shown in \cite{rus-sch}
that the problem of deciding whether a finite graph with loops allowed is homomorphism-homogeneous
is coNP-complete. In this section we adapt the argument of \cite{rus-sch} to show that
the same holds for bidirectionally connected improper digraphs. Consequently,
for the class of bidirectionally connected digraphs we cannot hope for a full description
that involves a catalogue.

\begin{figure}
  \centering
  \unitlength 0.24pt
\begin{picture}(300.0,200.0)
    \put(50.0,50.0){\circle*{12.0}}
    \put(250.0,50.0){\circle*{12.0}}
    \put(150.0,150.0){\circle*{12.0}}
    \path(50.0,50.0)(150.0,150.0)
    \path(150.0,150.0)(250.0,50.0)
    \put(150.0,168.0){\makebox(0,0)[b]{$\rightleftarrows$}}
    \put(50.0,32.0){\makebox(0,0)[t]{$\leftarrow$}}
    \put(250.0,32.0){\makebox(0,0)[t]{$\rightarrow$}}
\end{picture}
  \caption{The poset $\MM$}
  \label{nans.fig.poset}
\end{figure}
Let $M = \{\leftarrow, \rightarrow, \rightleftarrows\}$ and let
$\MM = (M, \le)$ be the three-element partially ordered
set depicted in Fig.~\ref{nans.fig.poset}.
Let $D$ be an improper digraph. We say that a vertex $c \in V(D)$ is a \emph{cone} for a sequence of
vertices $(u_1, \ldots, u_n) \in V(D)^n$ if $c \sim u_i$ for all $i$.
A vertex $c \in V(D)$ is a cone for the sequence of vertices $(u_1, \ldots, u_n)$
of \emph{type} $(t_1, \ldots, t_n) \in M^n$ if the following holds for every $i$:
\begin{itemize}
\item
  if $t_i =\mathstrut\rightarrow$ then $c \rightarrow u_i$,
\item
  if $t_i =\mathstrut\leftarrow$ then $u_i \rightarrow c$, and
\item
  if $t_i =\mathstrut\rightleftarrows$ then $c \rightleftarrows u_i$.
\end{itemize}
We say that a cone $c$ of type $(t_1, \ldots, t_n)$ for some sequence of vertices is \emph{not weaker than}
the cone $c'$ of type $(t'_1, \ldots, t'_n)$ for some (other) sequence of vertices
if $(t'_1, \ldots, t'_n) \le (t_1, \ldots, t_n)$. We write $c' \preccurlyeq c$.
The proof of the following lemma is analogous to the proof of \cite[Theorem 6]{rus-sch}:

\begin{LEM}\label{nans.lem.cones}
  A reflexive improper digraph $D$ is not homomorphism-homo\-gene\-ous if and only if there are vertices
  $u_1, \ldots, u_m, w_1, \ldots, w_m \in V(D)$ and a homomorphism $f : D[u_1, \ldots, u_m] \to D[w_1, \ldots, w_m]
  : u_i \mapsto w_i$ with the following property:
  \begin{itemize}
  \item
    either $(u_1, \ldots, u_m)$ has a cone and $(w_1, \ldots, w_m)$ does not,
  \item
    or $(u_1, \ldots, u_m)$ has a cone $c$ such that $c \not\preccurlyeq d$
    for every cone $d$ of $(w_1, \ldots, w_m)$.
  \end{itemize}
\end{LEM}

\begin{THM}\label{nans.thm.KK1}
  The problem of deciding whether an improper finite reflexive bidirectionally connected digraph is
  homomorphism-homogeneous is coNP-complete.
\end{THM}
\begin{PROOF}
  Let us first show that the problem is in coNP having in mind the criterion of homomorphism-homogeneity
  provided by Lemma~\ref{nans.lem.cones}.
  Given an improper digraph $D$ and a triple $((u_1, \ldots, u_m), (w_1, \ldots, w_m), f)$ where $f$ is
  a homomorphism from $D[u_1, \ldots, u_m]$ to $D[w_1, \ldots, w_m]$ that takes $u_i$ to $w_i$, one can check
  in polynomial time that
  \begin{itemize}
  \item
    either $(u_1, \ldots, u_m)$ has a cone and $(w_1, \ldots, w_m)$ does not,
  \item
    or $(u_1, \ldots, u_m)$ has a cone $c$ such that $c \not\preccurlyeq d$
    for every cone $d$ of $(w_1, \ldots, w_m)$.
  \end{itemize}

  We prove the hardness by reducing the INDEPENDENT SET problem.
  Take any integer $k \ge 2$ and an irreflexive graph $G = (V, E)$ where $V = \{v_1, \ldots, v_n\}$,
  and choose two $(k+1)$-element sets
  $I = \{q_0, q_1, \ldots, q_k\}$ and $S = \{s_0, s_1, \ldots, s_k\}$ in such a way that
  $V$, $I$ and $S$ are pairwise disjoint. Let $G_k$ be the reflexive improper digraph constructed as follows:
  \begin{align*}
    V(G_k) = &V \union I \union S\\
    E(G_k) = &E \union \{ v_i \to v_j : v_i \not\sim v_j \text{ in $G$ and } i < j\}
                 \union \{\text{all loops on } V \union I \union S\}\\
             &\union \{s_i \rightleftarrows s_j : i \ne j\} \union \{q_i \to q_j : i < j\}\\
             &\union \{ v \rightleftarrows q_j : v \in V, j > 0 \} \union \{v \to q_0 : v \in V \}\\
             &\union \{v \rightleftarrows s_j : v \in V, 0 \le j \le k\}\\
             &\union \{ s_i \rightleftarrows q_j : i \ne j\} \union \{s_i \to q_i : 0 \le i \le k \}.
  \end{align*}
  Note the following:
  \begin{itemize}
  \item
    $G$ is a graph, so if $v_i \sim v_j$ in $G$ then $v_i \rightleftarrows v_j$ in $G_k$;
  \item
    for every pair of distinct vertices $x, y \in V(G_k)$ we have either $x \to y$ or $y \to x$ or
    $x \rightleftarrows y$;
  \item
    if $X$ is an $m$-independent set in $G$, then $G_k[X]$ is a transitive tournament on $m$ vertices
    (modulo loops);
  \item
    $G_k[I]$ is a transitive tournament on $k+1$ vertices (modulo loops).
  \end{itemize}
  Let us show that $G$ has a $k$-independent set if and only if $G_k$ is not homo\-morph\-ism-homogeneous.

  $(\Rightarrow)$ Assume that $\{x_0, x_1, \ldots, x_{k-1}\} \subseteq V$ is a $k$-independent
  set in $G$. Then $G_k[x_0, x_1, \ldots, x_{k-1}]$ is a transitive tournament (modulo loops, of course)
  and without loss of generality we can assume that $x_0 \to x_1 \to \ldots \to x_{k-1}$.
  The mapping
  $$
    f = \begin{pmatrix}
      x_0 & x_1 & \ldots & x_{k-1} & q_0\\
      q_0 & q_1 & \ldots & q_{k-1} & q_k
    \end{pmatrix}.
  $$
  is a homomorphism from $G_k[x_0, x_1, \ldots, x_{k-1}, q_0]$ to $G_k[q_0, q_1, \ldots, q_k]$.
  If $G_k$ were homomorphism-homogeneous, then $f$ would extend to an endomorphism $f^*$ of $G_k$, so
  $f^*(s_1)$ would be a cone for $(q_0, q_1, \ldots, q_k)$ of type $(\rightleftarrows, \rightleftarrows, \ldots,
  \rightleftarrows)$ since $s_1$ is a cone for $(x_0, x_1, \ldots, x_{k-1}, q_0)$ of the same type.
  But it is easy to see that $(q_0, q_1, \ldots, q_k)$ does not have a cone of that type in~$G_k$.

  $(\Leftarrow)$ Assume that $G$ does not have a $k$-independent set and let us show that
  $G_k$ is homomorphism-homogeneous. Clearly, it suffices to show that every homomorphism $f : G_k[U] \to G_k[W]$
  where $W = f(U)$ can be extended to a homomorphism
  $f' : G_k[U \union \{x\}] \to G_k[W \union \{y\}]$ where $x \notin U$ and $f'(x) = y$.

  Take any homomorphism $f : G_k[U] \to G_k[W]$ where $W = f(U)$ and
  $U \ne V \union S \union I$. If $I \not\subseteq W$, say, $q_i \notin W$
  for some $q_i \in I$, then $s_i$ is a cone for $W$ in $G_k$ of type
  $(\rightleftarrows, \rightleftarrows, \ldots, \rightleftarrows)$. Now take any $x \notin U$ and note that
  $f' : G_k[U \union \{x\}] \to G_k[W \union \{s_i\}]$ where $f'(x) = s_i$ and $f'(y) = f(y)$ for $y \ne x$
  is a homomorphism which extends~$f$.

  Assume, now, that $I \subseteq W$. Then there exist $x_0, \ldots, x_k \in U$ such that
  $f(x_i) = q_i$, $0 \le i \le k$. Since $G_k[I]$ is a transitive tournament on $k+1$ vertices (modulo loops),
  so is $G_k[X]$ where $X = \{x_0, \ldots, x_k\}$.
  Let us show that $X \sec S = \0$. Clearly, $|X \sec S| \le 1$ since every pair of distinct vertices from $S$
  is connected by a double edge, while $G_k[X]$ is a tournament.
  If $|X \sec S| = 1$, say, $X \sec S = \{s_i\}$, then $X \sec V = \0$
  since every vertex from $S$ is connected by a double edge to every vertex from $V$. Therefore,
  $|X \sec I| = k \ge 2$, so there exists a $j \ne i$ such that $q_j \in X \sec I$.
  But, $q_j \rightleftarrows s_i$ by construction, which contradicts the fact that $G_k[X]$ is
  a tournament. This shows that $X \sec S = \0$.

  Next, let us show that $X \sec V = \0$. Assume this is not the case and let $v \in X \sec V$.
  Since $G$ does not have a $k$-independent set, it follows that no $k$-element subset of $V$ induces
  a tournament in $G_k$. So, $|X \sec V| \le k - 1$, whence $|X \sec I| \ge 2$. Consequently,
  there exists an $i > 0$ such that $q_i \in X$. But, $q_i \rightleftarrows v$ by construction,
  which contradicts the fact that $G_k[X]$ is a tournament. This shows that $X \sec V = \0$.

  Therefore, $X = I$ so $f(q_i) = q_i$, $0 \le i \le k$, since $G_k[I]$ is a transitive tournament.
  Moreover, the argument above shows that
  $$
  \text{if } f(x) \in I \text{ then } x \in I.
  \eqno{(\star)}
  $$

  If $V \not\subseteq U$, take any $v \in V \setminus U$
  and extend $f$ by setting $f'(v) = s_0$ and $f'(y) = f(y)$ for $y \in U$. Then $f'$ is
  a homomorphism (which clearly extends~$f$) since $v \to q_0$ and $s_{0} \to q_{0}$ by construction,
  while $s_{0} \rightleftarrows x$ for all $x \ne q_{0}$.
  If, however, $V \subseteq U$, then $S \not\subseteq U$. Take any $s_i \in S \setminus U$
  and extend $f$ by setting $f'(s_i) = s_{i}$ and $f'(y) = f(y)$ for $y \in U$. It is easy to see that
  $$
    f' = \Bigg(
      \overbrace{
        \begin{matrix}
          q_0 & \ldots & q_k \\
          q_0 & \ldots & q_k
        \end{matrix}
      }^I
      \;
      \overbrace{
        \begin{matrix}
          v_1 & \ldots & v_n \\
          z_1 & \ldots & z_n
        \end{matrix}
      }^V
      \;
      \overbrace{
        \begin{matrix}
          s_{j_1} & \ldots & s_{j_t} \\
          w_1 & \ldots & w_t
        \end{matrix}
      }^{U \sec S}
      \;
      \begin{matrix}
        s_i\\
        s_i
      \end{matrix}
    \Bigg)
  $$
  is a homomorphism from $G_k[U \union \{s_i\}]$ to $G_k[W \union \{s_{i}\}]$:
  $s_{i} \rightleftarrows x$ for all $x \ne q_{i}$, and $q_{i} \notin
  \{z_1, \ldots, z_n, w_1, \ldots, w_t\}$ because of~$(\star)$.
\end{PROOF}

\section{Digraphs with involution}
\label{nans.sec.dinv}

The classification of bidirectionally disconnected systems heavily relies on the following
peculiar class of digraphs.
Let $D$ be a reflexive improper digraph. We say that $D$ is a \emph{digraph with involution} if
there exists an automorphism $'$ of $D$ satisfying
\begin{itemize}
\item[(DI1)]
  $x = x''$;
\item[(DI2)]
  if $x \to y$ then $y \to x'$;
\item[(DI3)]
  if $x$ and $y$ are distinct vertices satisfying $x \rightleftarrows y$ then $y = x'$.
\end{itemize}

\begin{LEM}
  Let $D$ be a digraph with involution $'$.
  Then, for all $x, y \in V(D)$,

  $(a)$ $x \rightleftarrows x'$;

  $(b)$ $x = x'$ if and only if $x$ is an isolated vertex of $D$;

  $(c)$ if $x \ne y$ and $x \sim y$ then $\{x, x'\} \bowtie \{y, y'\}$.
\end{LEM}
\begin{PROOF}
  $(a)$ Since $D$ is reflexive, we have $x \to x$ and $x' \to x'$. From (DI2) we now conclude
  $x \to x'$ and $x' \to x'' = x$.

  $(b)$ Assume that $x$ is an isolated vertex of $D$. Then by $(a)$ we have $x \rightleftarrows x'$
  whence $x = x'$. For the converse, assume that $x = x'$ and let us show that $x$ is then an isolated
  vertex of $D$. Suppose this is not the case, and let $y$ be a vertex distinct from $x$ such that
  $x \sim y$, say $x \to y$. Then by (DI2) we conclude that $y \to x' = x$. Therefore, $x \rightleftarrows y$,
  so (DI3) now yields $y = x' = x$, which contradicts the assumption $y \ne x$.

  $(c)$ Assume that $x \to y$. Then $y \to x'$ by (DI2), $x' \to y'$ since $'$ is an automorphism of $D$
  and, by the same argument, $y' \to x'' = x$.
\end{PROOF}

Clearly, if $D$ is a digraph with involution $'$ then each class $S$ of $\theta(D)$ takes the
form $\{x, x'\}$ (see Fig.~\ref{nans.fig.skeleton}~$(a)$). So we have the following:
\begin{figure}
  \centering
  \unitlength 0.24pt
\begin{picture}(1425.0,900.0)
    \put(225.0,200.0){\circle*{12.0}}
    \put(75.0,300.0){\circle*{12.0}}
    \put(250.0,400.0){\arc{403.113}{1.69513}{2.62243}}
    \put(50.0,100.0){\arc{403.113}{4.83672}{5.76402}}
    \put(127.398,243.281){\arc{77.4144}{5.20764}{5.81636}}
    \put(164.191,311.393){\arc{77.4144}{1.45733}{2.06605}}
    \put(172.602,256.719){\arc{77.4144}{2.06605}{2.67477}}
    \put(135.809,188.607){\arc{77.4144}{4.59892}{5.20764}}
    \put(135.809,388.607){\arc{77.4144}{4.59892}{5.20764}}
    \put(172.602,456.719){\arc{77.4144}{2.06605}{2.67477}}
    \put(164.191,511.393){\arc{77.4144}{1.45733}{2.06605}}
    \put(127.398,443.281){\arc{77.4144}{5.20764}{5.81636}}
    \put(50.0,300.0){\arc{403.113}{4.83672}{5.76402}}
    \put(250.0,600.0){\arc{403.113}{1.69513}{2.62243}}
    \put(75.0,500.0){\circle*{12.0}}
    \put(225.0,400.0){\circle*{12.0}}
    \put(135.809,688.607){\arc{77.4144}{4.59892}{5.20764}}
    \put(172.602,756.719){\arc{77.4144}{2.06605}{2.67477}}
    \put(164.191,811.393){\arc{77.4144}{1.45733}{2.06605}}
    \put(127.398,743.281){\arc{77.4144}{5.20764}{5.81636}}
    \put(50.0,600.0){\arc{403.113}{4.83672}{5.76402}}
    \put(250.0,900.0){\arc{403.113}{1.69513}{2.62243}}
    \put(75.0,800.0){\circle*{12.0}}
    \put(225.0,700.0){\circle*{12.0}}
    \put(385.809,488.607){\arc{77.4144}{4.59892}{5.20764}}
    \put(422.602,556.719){\arc{77.4144}{2.06605}{2.67477}}
    \put(414.191,611.393){\arc{77.4144}{1.45733}{2.06605}}
    \put(377.398,543.281){\arc{77.4144}{5.20764}{5.81636}}
    \put(300.0,400.0){\arc{403.113}{4.83672}{5.76402}}
    \put(500.0,700.0){\arc{403.113}{1.69513}{2.62243}}
    \put(325.0,600.0){\circle*{12.0}}
    \put(475.0,500.0){\circle*{12.0}}
    \path(75.0,300.0)(75.0,500.0)
    \path(75.0,500.0)(225.0,200.0)
    \path(225.0,200.0)(225.0,400.0)
    \path(225.0,400.0)(75.0,300.0)
    \path(75.0,500.0)(75.0,800.0)
    \path(75.0,800.0)(225.0,400.0)
    \path(225.0,400.0)(225.0,700.0)
    \path(225.0,700.0)(75.0,500.0)
    \path(325.0,600.0)(75.0,500.0)
    \path(75.0,500.0)(475.0,500.0)
    \path(475.0,500.0)(225.0,400.0)
    \path(225.0,400.0)(325.0,600.0)
    \put(189.962,322.558){\arc{90.0}{4.12437}{4.64797}}
    \put(140.038,397.442){\arc{90.0}{0.459175}{0.982774}}
    \put(180.0,280.0){\arc{90.0}{-2.00012e-005}{0.523579}}
    \put(270.0,280.0){\arc{90.0}{2.61797}{3.14157}}
    \put(175.249,400.125){\arc{90.0}{2.67793}{3.20152}}
    \put(94.7508,359.875){\arc{90.0}{5.29592}{5.81952}}
    \put(30.0,380.0){\arc{90.0}{-2.00012e-005}{0.523579}}
    \put(120.0,380.0){\arc{90.0}{2.61797}{3.14157}}
    \put(177.135,655.801){\arc{90.0}{2.7828}{3.3064}}
    \put(92.8652,624.199){\arc{90.0}{5.4008}{5.92439}}
    \put(201.0,593.0){\arc{90.0}{3.78507}{4.30867}}
    \put(129.0,647.0){\arc{90.0}{0.119882}{0.643481}}
    \put(180.0,635.732){\arc{90.0}{-2.00012e-005}{0.523579}}
    \put(270.0,635.732){\arc{90.0}{2.61797}{3.14157}}
    \put(30.0,680.0){\arc{90.0}{-2.00012e-005}{0.523579}}
    \put(120.0,680.0){\arc{90.0}{2.61797}{3.14157}}
    \put(315.0,545.0){\arc{90.0}{1.57078}{2.09438}}
    \put(315.0,455.0){\arc{90.0}{4.18877}{4.71237}}
    \put(262.051,574.725){\arc{90.0}{0.463628}{0.987226}}
    \put(342.55,534.476){\arc{90.0}{3.08162}{3.60522}}
    \put(366.713,408.219){\arc{90.0}{4.33186}{4.85546}}
    \put(333.287,491.781){\arc{90.0}{0.666671}{1.19027}}
    \put(151.505,482.135){\arc{90.0}{4.33186}{4.85546}}
    \put(118.08,565.698){\arc{90.0}{0.666671}{1.19027}}
    \put(435.809,176.107){\arc{77.4144}{4.59892}{5.20764}}
    \put(472.602,244.219){\arc{77.4144}{2.06605}{2.67477}}
    \put(464.191,298.893){\arc{77.4144}{1.45733}{2.06605}}
    \put(427.398,230.781){\arc{77.4144}{5.20764}{5.81636}}
    \put(350.0,87.5){\arc{403.113}{4.83672}{5.76402}}
    \put(550.0,387.5){\arc{403.113}{1.69513}{2.62243}}
    \put(375.0,287.5){\circle*{12.0}}
    \put(525.0,187.5){\circle*{12.0}}
    \put(350.0,87.5){\circle*{12.0}}
    \put(225.0,38.0){\makebox(0,0)[t]{$(a)$}}
    \put(1075.0,38.0){\makebox(0,0)[t]{$(b)$}}
    \put(800.0,200.0){\circle*{12.0}}
    \put(900.0,100.0){\circle*{12.0}}
    \put(943.75,243.75){\arc{300.52}{1.86622}{2.84613}}
    \put(756.25,56.25){\arc{300.52}{5.00781}{5.98772}}
    \put(890.146,143.597){\arc{73.3281}{2.0622}{2.70484}}
    \put(855.544,78.946){\arc{73.3281}{4.56115}{5.20379}}
    \put(1294.46,671.054){\arc{73.3281}{1.41956}{2.0622}}
    \put(1259.85,606.403){\arc{73.3281}{5.20379}{5.84644}}
    \put(1206.25,506.25){\arc{300.52}{5.00781}{5.98772}}
    \put(1393.75,693.75){\arc{300.52}{1.86622}{2.84613}}
    \put(1350.0,550.0){\circle*{12.0}}
    \put(1250.0,650.0){\circle*{12.0}}
    \put(900.0,650.0){\circle*{12.0}}
    \put(800.0,550.0){\circle*{12.0}}
    \put(756.25,693.75){\arc{300.52}{0.295421}{1.27534}}
    \put(943.75,506.25){\arc{300.52}{3.43701}{4.41693}}
    \put(890.146,606.403){\arc{73.3281}{3.5783}{4.22095}}
    \put(855.544,671.054){\arc{73.3281}{1.07935}{1.722}}
    \put(1294.46,78.946){\arc{73.3281}{4.22095}{4.86359}}
    \put(1259.85,143.597){\arc{73.3281}{0.436708}{1.07935}}
    \put(1393.75,56.25){\arc{300.52}{3.43701}{4.41693}}
    \put(1206.25,243.75){\arc{300.52}{0.295421}{1.27534}}
    \put(1250.0,100.0){\circle*{12.0}}
    \put(1350.0,200.0){\circle*{12.0}}
    \path(900.0,650.0)(1250.0,650.0)
    \path(800.0,550.0)(1350.0,550.0)
    \path(1250.0,650.0)(800.0,550.0)
    \path(1350.0,550.0)(900.0,650.0)
    \path(900.0,100.0)(1250.0,100.0)
    \path(1250.0,100.0)(800.0,200.0)
    \path(800.0,200.0)(1350.0,200.0)
    \path(1350.0,200.0)(900.0,100.0)
    \path(800.0,200.0)(800.0,550.0)
    \path(800.0,550.0)(900.0,100.0)
    \path(900.0,100.0)(900.0,650.0)
    \path(900.0,650.0)(800.0,200.0)
    \path(1350.0,200.0)(1350.0,550.0)
    \path(1350.0,550.0)(1250.0,100.0)
    \path(1250.0,100.0)(1250.0,650.0)
    \path(1250.0,650.0)(1350.0,200.0)
    \path(800.0,200.0)(1250.0,650.0)
    \path(1250.0,650.0)(900.0,100.0)
    \path(900.0,100.0)(1350.0,550.0)
    \path(1350.0,550.0)(800.0,200.0)
    \path(900.0,650.0)(1350.0,200.0)
    \path(1350.0,200.0)(800.0,550.0)
    \path(800.0,550.0)(1250.0,100.0)
    \path(1250.0,100.0)(900.0,650.0)
    \put(980.238,586.072){\arc{90.0}{4.93104}{5.45464}}
    \put(999.762,673.928){\arc{90.0}{1.26585}{1.78945}}
    \put(1180.0,695.0){\arc{90.0}{1.57078}{2.09438}}
    \put(1180.0,605.0){\arc{90.0}{4.18877}{4.71237}}
    \put(989.762,546.072){\arc{90.0}{4.4937}{5.0173}}
    \put(970.238,633.928){\arc{90.0}{0.828509}{1.35211}}
    \put(1130.0,595.0){\arc{90.0}{1.57078}{2.09438}}
    \put(1130.0,505.0){\arc{90.0}{4.18877}{4.71237}}
    \put(755.0,480.0){\arc{90.0}{0.0}{0.523579}}
    \put(845.0,480.0){\arc{90.0}{2.61797}{3.14157}}
    \put(1305.0,480.0){\arc{90.0}{0.0}{0.523579}}
    \put(1395.0,480.0){\arc{90.0}{2.61797}{3.14157}}
    \put(1180.0,145.0){\arc{90.0}{1.57078}{2.09438}}
    \put(1180.0,55.0){\arc{90.0}{4.18877}{4.71237}}
    \put(1373.93,299.762){\arc{90.0}{2.9229}{3.4465}}
    \put(1286.07,280.238){\arc{90.0}{5.5409}{6.0645}}
    \put(863.928,280.238){\arc{90.0}{3.36024}{3.88384}}
    \put(776.072,299.762){\arc{90.0}{5.97824}{6.50181}}
    \put(999.762,76.0716){\arc{90.0}{4.4937}{5.0173}}
    \put(980.238,163.928){\arc{90.0}{0.828509}{1.35211}}
    \put(855.0,430.0){\arc{90.0}{0.0}{0.523579}}
    \put(945.0,430.0){\arc{90.0}{2.61797}{3.14157}}
    \put(903.928,289.762){\arc{90.0}{2.9229}{3.4465}}
    \put(816.072,270.238){\arc{90.0}{5.5409}{6.0645}}
    \put(1130.0,245.0){\arc{90.0}{1.57078}{2.09438}}
    \put(1130.0,155.0){\arc{90.0}{4.18877}{4.71237}}
    \put(970.238,116.072){\arc{90.0}{4.93104}{5.45464}}
    \put(989.762,203.928){\arc{90.0}{1.26585}{1.78945}}
    \put(1205.0,430.0){\arc{90.0}{0.0}{0.523579}}
    \put(1295.0,430.0){\arc{90.0}{2.61797}{3.14157}}
    \put(1333.93,270.238){\arc{90.0}{3.36024}{3.88384}}
    \put(1246.07,289.762){\arc{90.0}{5.97824}{6.50181}}
    \put(1083.18,546.82){\arc{90.0}{0.785378}{1.30898}}
    \put(1146.82,483.18){\arc{90.0}{3.40337}{3.92697}}
    \put(1183.18,446.82){\arc{90.0}{0.785378}{1.30898}}
    \put(1246.82,383.18){\arc{90.0}{3.40337}{3.92697}}
    \put(1246.82,366.82){\arc{90.0}{2.35617}{2.87977}}
    \put(1183.18,303.18){\arc{90.0}{4.97417}{5.49777}}
    \put(1146.82,266.82){\arc{90.0}{2.35617}{2.87977}}
    \put(1083.18,203.18){\arc{90.0}{4.97417}{5.49777}}
    \put(1007.96,185.841){\arc{90.0}{3.7083}{4.2319}}
    \put(932.035,234.159){\arc{90.0}{0.0431104}{0.566709}}
    \put(934.159,232.035){\arc{90.0}{4.14564}{4.66924}}
    \put(885.841,307.965){\arc{90.0}{0.480448}{1.00405}}
    \put(885.841,442.035){\arc{90.0}{5.2791}{5.8027}}
    \put(934.159,517.965){\arc{90.0}{1.61391}{2.13751}}
    \put(932.035,515.841){\arc{90.0}{5.71644}{6.24003}}
    \put(1007.96,564.159){\arc{90.0}{2.05124}{2.57484}}
    \put(828.331,144.048){\arc{73.3281}{5.39978}{6.04242}}
    \put(874.86,200.723){\arc{73.3281}{1.61554}{2.25818}}
    \put(1321.67,605.952){\arc{73.3281}{2.25818}{2.90083}}
    \put(1275.14,549.277){\arc{73.3281}{4.75713}{5.39978}}
    \put(1275.14,200.723){\arc{73.3281}{0.88337}{1.52601}}
    \put(1321.67,144.048){\arc{73.3281}{3.38232}{4.02496}}
    \put(874.86,549.277){\arc{73.3281}{4.02496}{4.66761}}
    \put(828.331,605.952){\arc{73.3281}{0.240725}{0.88337}}
    \path(96.6506,837.5)(75.0,800.0)
    \put(75.0,850.0){\arc{50.0}{2.61797}{6.80674}}
    \path(75.0,800.0)(53.3494,837.5)
    \put(44.6891,837.5){\arc{90.0}{0.0}{0.523579}}
    \put(122.631,792.5){\arc{90.0}{3.66517}{4.18877}}
    \path(346.651,637.5)(325.0,600.0)
    \put(325.0,650.0){\arc{50.0}{2.61797}{6.80674}}
    \path(325.0,600.0)(303.349,637.5)
    \put(294.689,637.5){\arc{90.0}{0.0}{0.523579}}
    \put(372.631,592.5){\arc{90.0}{3.66517}{4.18877}}
    \path(396.651,325.0)(375.0,287.5)
    \put(375.0,337.5){\arc{50.0}{2.61797}{6.80674}}
    \path(375.0,287.5)(353.349,325.0)
    \put(344.689,325.0){\arc{90.0}{0.0}{0.523579}}
    \put(422.631,280.0){\arc{90.0}{3.66517}{4.18877}}
    \path(37.5,521.651)(75.0,500.0)
    \put(25.0,500.0){\arc{50.0}{1.04718}{5.23597}}
    \path(75.0,500.0)(37.5,478.349)
    \put(37.5,469.689){\arc{90.0}{4.71237}{5.23597}}
    \put(82.5,547.631){\arc{90.0}{2.09438}{2.61797}}
    \path(37.5,321.651)(75.0,300.0)
    \put(25.0,300.0){\arc{50.0}{1.04718}{5.23597}}
    \path(75.0,300.0)(37.5,278.349)
    \put(37.5,269.689){\arc{90.0}{4.71237}{5.23597}}
    \put(82.5,347.631){\arc{90.0}{2.09438}{2.61797}}
    \path(203.349,162.5)(225.0,200.0)
    \put(225.0,150.0){\arc{50.0}{5.75957}{9.94834}}
    \path(225.0,200.0)(246.651,162.5)
    \put(255.311,162.5){\arc{90.0}{3.14157}{3.66517}}
    \put(177.369,207.5){\arc{90.0}{0.523579}{1.04718}}
    \path(503.349,150.0)(525.0,187.5)
    \put(525.0,137.5){\arc{50.0}{5.75957}{9.94834}}
    \path(525.0,187.5)(546.651,150.0)
    \put(555.311,150.0){\arc{90.0}{3.14157}{3.66517}}
    \put(477.369,195.0){\arc{90.0}{0.523579}{1.04718}}
    \path(453.349,462.5)(475.0,500.0)
    \put(475.0,450.0){\arc{50.0}{5.75957}{9.94834}}
    \path(475.0,500.0)(496.651,462.5)
    \put(505.311,462.5){\arc{90.0}{3.14157}{3.66517}}
    \put(427.369,507.5){\arc{90.0}{0.523579}{1.04718}}
    \path(266.826,711.207)(225.0,700.0)
    \put(260.355,735.355){\arc{50.0}{3.40337}{7.59214}}
    \put(230.083,747.95){\arc{90.0}{0.785378}{1.30898}}
    \put(253.377,661.016){\arc{90.0}{4.45057}{4.97417}}
    \path(225.0,700.0)(236.207,741.826)
    \path(239.81,359.31)(225.0,400.0)
    \put(263.302,367.861){\arc{50.0}{4.8869}{9.07567}}
    \path(225.0,400.0)(267.643,392.481)
    \put(273.21,399.115){\arc{90.0}{2.26891}{2.79251}}
    \put(188.638,368.333){\arc{90.0}{5.9341}{6.45768}}
    \put(397.631,80.0){\arc{90.0}{3.66517}{4.18877}}
    \put(319.689,125.0){\arc{90.0}{0.0}{0.523579}}
    \path(350.0,87.5)(328.349,125.0)
    \put(350.0,137.5){\arc{50.0}{2.61797}{6.80674}}
    \path(371.651,125.0)(350.0,87.5)
    \put(947.631,642.5){\arc{90.0}{3.66517}{4.18877}}
    \put(869.689,687.5){\arc{90.0}{0.0}{0.523579}}
    \path(900.0,650.0)(878.349,687.5)
    \put(900.0,700.0){\arc{50.0}{2.61797}{6.80674}}
    \path(921.651,687.5)(900.0,650.0)
    \put(1297.63,642.5){\arc{90.0}{3.66517}{4.18877}}
    \put(1219.69,687.5){\arc{90.0}{0.0}{0.523579}}
    \path(1250.0,650.0)(1228.35,687.5)
    \put(1250.0,700.0){\arc{50.0}{2.61797}{6.80674}}
    \path(1271.65,687.5)(1250.0,650.0)
    \put(807.5,597.631){\arc{90.0}{2.09438}{2.61797}}
    \put(762.5,519.689){\arc{90.0}{4.71237}{5.23597}}
    \path(800.0,550.0)(762.5,528.349)
    \put(750.0,550.0){\arc{50.0}{1.04718}{5.23597}}
    \path(762.5,571.651)(800.0,550.0)
    \put(807.5,247.631){\arc{90.0}{2.09438}{2.61797}}
    \put(762.5,169.689){\arc{90.0}{4.71237}{5.23597}}
    \path(800.0,200.0)(762.5,178.349)
    \put(750.0,200.0){\arc{50.0}{1.04718}{5.23597}}
    \path(762.5,221.651)(800.0,200.0)
    \put(852.369,107.5){\arc{90.0}{0.523579}{1.04718}}
    \put(930.311,62.5){\arc{90.0}{3.14157}{3.66517}}
    \path(900.0,100.0)(921.651,62.5)
    \put(900.0,50.0){\arc{50.0}{5.75957}{9.94834}}
    \path(878.349,62.5)(900.0,100.0)
    \put(1202.37,107.5){\arc{90.0}{0.523579}{1.04718}}
    \put(1280.31,62.5){\arc{90.0}{3.14157}{3.66517}}
    \path(1250.0,100.0)(1271.65,62.5)
    \put(1250.0,50.0){\arc{50.0}{5.75957}{9.94834}}
    \path(1228.35,62.5)(1250.0,100.0)
    \put(1342.5,152.369){\arc{90.0}{5.23597}{5.75957}}
    \put(1387.5,230.311){\arc{90.0}{1.57078}{2.09438}}
    \path(1350.0,200.0)(1387.5,221.651)
    \put(1400.0,200.0){\arc{50.0}{4.18877}{8.37754}}
    \path(1387.5,178.349)(1350.0,200.0)
    \put(1342.5,502.369){\arc{90.0}{5.23597}{5.75957}}
    \put(1387.5,580.311){\arc{90.0}{1.57078}{2.09438}}
    \path(1350.0,550.0)(1387.5,571.651)
    \put(1400.0,550.0){\arc{50.0}{4.18877}{8.37754}}
    \path(1387.5,528.349)(1350.0,550.0)
\end{picture}
  \caption{$(a)$ A digraph with involution; $(b)$ A tournament with involution}
  \label{nans.fig.skeleton}
\end{figure}
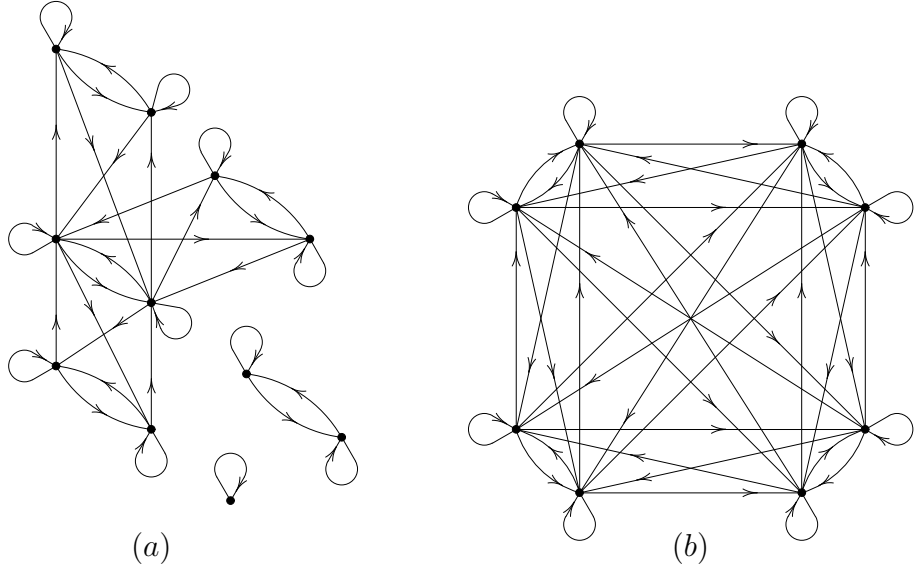

\begin{COR}
  If $D$ is a digraph with involution, then the automorphism $'$ of $D$ satisfying
  (DI1), (DI2) and (DI3) is unique.
\end{COR}

A digraph with involution is \emph{a tournament with involution} if $x \sim y$
for all $x, y \in V(D)$ (Fig.~\ref{nans.fig.skeleton}~$(b)$).

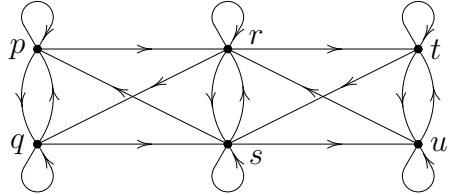
\begin{figure}
  \centering
  \unitlength 0.24pt
\begin{picture}(775.0,300.0)
    \put(75.0,225.0){\circle*{12.0}}
    \put(75.0,75.0){\circle*{12.0}}
    \put(175.0,150.0){\arc{250.0}{2.49807}{3.78507}}
    \put(-25.0,150.0){\arc{250.0}{5.63966}{6.92665}}
    \put(275.0,150.0){\arc{250.0}{5.63966}{6.92665}}
    \put(475.0,150.0){\arc{250.0}{2.49807}{3.78507}}
    \put(375.0,75.0){\circle*{12.0}}
    \put(375.0,225.0){\circle*{12.0}}
    \put(385.844,138.462){\arc{70.2011}{3.01287}{3.68414}}
    \put(316.224,129.452){\arc{70.2011}{5.4832}{6.15447}}
    \put(364.156,161.537){\arc{70.2011}{6.15448}{6.82573}}
    \put(433.777,170.546){\arc{70.2011}{2.34161}{3.01288}}
    \put(64.1561,161.537){\arc{70.2011}{6.15448}{6.82573}}
    \put(133.777,170.546){\arc{70.2011}{2.34161}{3.01288}}
    \put(85.8441,138.462){\arc{70.2011}{3.01287}{3.68414}}
    \put(16.2235,129.452){\arc{70.2011}{5.4832}{6.15447}}
    \path(75.0,225.0)(375.0,225.0)
    \path(375.0,225.0)(75.0,75.0)
    \path(75.0,75.0)(375.0,75.0)
    \path(375.0,75.0)(75.0,225.0)
    \put(255.0,120.0){\arc{90.0}{1.57078}{2.09438}}
    \put(255.0,30.0){\arc{90.0}{4.18877}{4.71237}}
    \put(255.0,270.0){\arc{90.0}{1.57078}{2.09438}}
    \put(255.0,180.0){\arc{90.0}{4.18877}{4.71237}}
    \put(174.875,124.751){\arc{90.0}{5.17602}{5.69962}}
    \put(215.125,205.249){\arc{90.0}{1.51083}{2.03442}}
    \put(275.125,124.751){\arc{90.0}{4.24872}{4.77232}}
    \put(234.875,205.249){\arc{90.0}{0.58353}{1.10713}}
    \put(534.875,205.249){\arc{90.0}{0.58353}{1.10713}}
    \put(575.125,124.751){\arc{90.0}{4.24872}{4.77232}}
    \put(515.125,205.249){\arc{90.0}{1.51083}{2.03442}}
    \put(474.875,124.751){\arc{90.0}{5.17602}{5.69962}}
    \put(555.0,180.0){\arc{90.0}{4.18877}{4.71237}}
    \put(555.0,270.0){\arc{90.0}{1.57078}{2.09438}}
    \put(555.0,30.0){\arc{90.0}{4.18877}{4.71237}}
    \put(555.0,120.0){\arc{90.0}{1.57078}{2.09438}}
    \path(675.0,75.0)(375.0,225.0)
    \path(375.0,75.0)(675.0,75.0)
    \path(675.0,225.0)(375.0,75.0)
    \path(375.0,225.0)(675.0,225.0)
    \put(733.777,170.546){\arc{70.2011}{2.34161}{3.01288}}
    \put(664.156,161.537){\arc{70.2011}{6.15448}{6.82573}}
    \put(616.224,129.452){\arc{70.2011}{5.4832}{6.15447}}
    \put(685.844,138.462){\arc{70.2011}{3.01287}{3.68414}}
    \put(675.0,225.0){\circle*{12.0}}
    \put(675.0,75.0){\circle*{12.0}}
    \put(775.0,150.0){\arc{250.0}{2.49807}{3.78507}}
    \put(575.0,150.0){\arc{250.0}{5.63966}{6.92665}}
    \put(122.631,217.5){\arc{90.0}{3.66517}{4.18877}}
    \put(44.6891,262.5){\arc{90.0}{0.0}{0.523579}}
    \path(75.0,225.0)(53.3494,262.5)
    \put(75.0,275.0){\arc{50.0}{2.61797}{6.80674}}
    \path(96.6506,262.5)(75.0,225.0)
    \put(422.631,217.5){\arc{90.0}{3.66517}{4.18877}}
    \put(344.689,262.5){\arc{90.0}{0.0}{0.523579}}
    \path(375.0,225.0)(353.349,262.5)
    \put(375.0,275.0){\arc{50.0}{2.61797}{6.80674}}
    \path(396.651,262.5)(375.0,225.0)
    \put(722.631,217.5){\arc{90.0}{3.66517}{4.18877}}
    \put(644.689,262.5){\arc{90.0}{0.0}{0.523579}}
    \path(675.0,225.0)(653.349,262.5)
    \put(675.0,275.0){\arc{50.0}{2.61797}{6.80674}}
    \path(696.651,262.5)(675.0,225.0)
    \path(696.651,37.5)(675.0,75.0)
    \put(675.0,25.0){\arc{50.0}{5.75957}{9.94834}}
    \path(675.0,75.0)(653.349,37.5)
    \put(644.689,37.5){\arc{90.0}{5.75957}{6.28317}}
    \put(722.631,82.5){\arc{90.0}{2.09438}{2.61797}}
    \path(396.651,37.5)(375.0,75.0)
    \put(375.0,25.0){\arc{50.0}{5.75957}{9.94834}}
    \path(375.0,75.0)(353.349,37.5)
    \put(344.689,37.5){\arc{90.0}{5.75957}{6.28317}}
    \put(422.631,82.5){\arc{90.0}{2.09438}{2.61797}}
    \path(96.6506,37.5)(75.0,75.0)
    \put(75.0,25.0){\arc{50.0}{5.75957}{9.94834}}
    \path(75.0,75.0)(53.3494,37.5)
    \put(44.6891,37.5){\arc{90.0}{5.75957}{6.28317}}
    \put(122.631,82.5){\arc{90.0}{2.09438}{2.61797}}
    \put(55.8771,225.0){\makebox(0,0)[r]{$p$}}
    \put(55.8771,75.0){\makebox(0,0)[r]{$q$}}
    \put(694.123,225.0){\makebox(0,0)[l]{$t$}}
    \put(694.123,75.0){\makebox(0,0)[l]{$u$}}
    \put(408.0,233.0){\makebox(0,0)[bl]{$r$}}
    \put(408.0,65.4918){\makebox(0,0)[tl]{$s$}}
\end{picture}
  \caption{The proof of Lemma \ref{nans.lem.dprime-hh}}
  \label{nans.fig.infl-skel}
\end{figure}
\begin{LEM}\label{nans.lem.dprime-hh}
  Let $D$ be a homomorphism-homogeneous digraph with involution.
  Then, for every connected component $S$ of $D$, we have that $D[S]$ is a tournament
  with involution.
\end{LEM}
\begin{PROOF}
  Suppose that there exists a connected component of $D$ such that $D[S]$ is
  not a tournament with involution.
  Then there exist $\{p, q\}, \{r, s\}, \{t, u\} \in V(D) / \theta(D)$
  such that $\{p, q\} \bowtie \{r, s\} \bowtie \{t, u\}$, but $\lnot(\{p, q\} \bowtie \{t, u\})$.
  Hence, $\{p, q\} \not\sim \{t, u\}$.
  Without loss of generality we can assume that $p \to r \to q \to s \to p$ and
  $r \to t \to s \to u \to r$, Fig.~\ref{nans.fig.infl-skel}. The mapping
  $$
    f : \begin{pmatrix}
      q & t & u \\
      q & t & t
    \end{pmatrix}
  $$
  is a homomorphism between finitely induced substructures of $D$ so it extends to an endomorphism $f^*$ of $D$.
  From $u \to r \to t$ it follows that $f^*(r) \rightleftarrows t$, so $f^*(r) \in \{u, t\}$.
  On the other hand, $r \to q$ implies $f^*(r) \to q$. Therefore, $\{t, u\} \sim \{p, q\}$. Contradiction.
\end{PROOF}

\begin{figure}
  \centering
  \unitlength 0.24pt
\begin{picture}(1500.0,400.0)
    \path(50.0,150.0)(250.0,150.0)
    \path(400.0,150.0)(600.0,150.0)
    \path(600.0,150.0)(500.0,300.0)
    \path(400.0,150.0)(500.0,300.0)
    \put(469.689,337.5){\arc{90.0}{0.0}{0.523579}}
    \path(500.0,300.0)(478.349,337.5)
    \path(521.651,337.5)(500.0,300.0)
    \put(500.0,350.0){\arc{50.0}{2.61797}{6.80674}}
    \put(547.631,292.5){\arc{90.0}{3.66517}{4.18877}}
    \put(19.6891,187.5){\arc{90.0}{0.0}{0.523579}}
    \path(50.0,150.0)(28.3494,187.5)
    \path(71.6506,187.5)(50.0,150.0)
    \put(50.0,200.0){\arc{50.0}{2.61797}{6.80674}}
    \put(97.6314,142.5){\arc{90.0}{3.66517}{4.18877}}
    \put(219.689,187.5){\arc{90.0}{0.0}{0.523579}}
    \path(250.0,150.0)(228.349,187.5)
    \path(271.651,187.5)(250.0,150.0)
    \put(250.0,200.0){\arc{50.0}{2.61797}{6.80674}}
    \put(297.631,142.5){\arc{90.0}{3.66517}{4.18877}}
    \put(422.558,264.962){\arc{90.0}{0.587983}{1.11158}}
    \put(497.442,215.038){\arc{90.0}{3.20598}{3.72958}}
    \put(502.558,215.038){\arc{90.0}{5.69516}{6.21876}}
    \put(577.442,264.962){\arc{90.0}{2.02997}{2.55357}}
    \put(520.0,195.0){\arc{90.0}{1.57078}{2.09438}}
    \put(520.0,105.0){\arc{90.0}{4.18877}{4.71237}}
    \put(170.0,195.0){\arc{90.0}{1.57078}{2.09438}}
    \put(170.0,105.0){\arc{90.0}{4.18877}{4.71237}}
    \put(50.0,150.0){\circle*{12.0}}
    \put(250.0,150.0){\circle*{12.0}}
    \put(400.0,150.0){\circle*{12.0}}
    \put(600.0,150.0){\circle*{12.0}}
    \put(500.0,300.0){\circle*{12.0}}
    \put(447.631,157.5){\arc{90.0}{2.09438}{2.61797}}
    \put(400.0,100.0){\arc{50.0}{5.75957}{9.94834}}
    \path(421.651,112.5)(400.0,150.0)
    \path(400.0,150.0)(378.349,112.5)
    \put(369.689,112.5){\arc{90.0}{5.75957}{6.28317}}
    \put(647.631,157.5){\arc{90.0}{2.09438}{2.61797}}
    \put(600.0,100.0){\arc{50.0}{5.75957}{9.94834}}
    \path(621.651,112.5)(600.0,150.0)
    \path(600.0,150.0)(578.349,112.5)
    \put(569.689,112.5){\arc{90.0}{5.75957}{6.28317}}
    \put(919.689,112.5){\arc{90.0}{5.75957}{6.28317}}
    \path(950.0,150.0)(928.349,112.5)
    \path(971.651,112.5)(950.0,150.0)
    \put(950.0,100.0){\arc{50.0}{5.75957}{9.94834}}
    \put(997.631,157.5){\arc{90.0}{2.09438}{2.61797}}
    \put(719.689,112.5){\arc{90.0}{5.75957}{6.28317}}
    \path(750.0,150.0)(728.349,112.5)
    \path(771.651,112.5)(750.0,150.0)
    \put(750.0,100.0){\arc{50.0}{5.75957}{9.94834}}
    \put(797.631,157.5){\arc{90.0}{2.09438}{2.61797}}
    \put(850.0,300.0){\circle*{12.0}}
    \put(950.0,150.0){\circle*{12.0}}
    \put(750.0,150.0){\circle*{12.0}}
    \put(870.0,105.0){\arc{90.0}{4.18877}{4.71237}}
    \put(870.0,195.0){\arc{90.0}{1.57078}{2.09438}}
    \put(904.016,300.101){\arc{90.0}{2.02997}{2.55357}}
    \put(829.132,250.178){\arc{90.0}{5.69516}{6.21876}}
    \put(847.442,215.038){\arc{90.0}{3.20598}{3.72958}}
    \put(772.558,264.962){\arc{90.0}{0.587983}{1.11158}}
    \put(897.631,292.5){\arc{90.0}{3.66517}{4.18877}}
    \put(850.0,350.0){\arc{50.0}{2.61797}{6.80674}}
    \path(871.651,337.5)(850.0,300.0)
    \path(850.0,300.0)(828.349,337.5)
    \put(819.689,337.5){\arc{90.0}{0.0}{0.523579}}
    \path(750.0,150.0)(850.0,300.0)
    \path(950.0,150.0)(850.0,300.0)
    \path(750.0,150.0)(950.0,150.0)
    \path(850.0,300.0)(1050.0,300.0)
    \path(950.0,150.0)(1050.0,300.0)
    \path(750.0,150.0)(1050.0,300.0)
    \put(1019.69,337.5){\arc{90.0}{0.0}{0.523579}}
    \path(1050.0,300.0)(1028.35,337.5)
    \path(1071.65,337.5)(1050.0,300.0)
    \put(1050.0,350.0){\arc{50.0}{2.61797}{6.80674}}
    \put(1097.63,292.5){\arc{90.0}{3.66517}{4.18877}}
    \put(1050.0,300.0){\circle*{12.0}}
    \put(972.558,264.962){\arc{90.0}{0.587983}{1.11158}}
    \put(1047.44,215.038){\arc{90.0}{3.20598}{3.72958}}
    \put(970.0,345.0){\arc{90.0}{1.57078}{2.09438}}
    \put(970.0,255.0){\arc{90.0}{4.18877}{4.71237}}
    \put(969.875,310.249){\arc{90.0}{1.10713}{1.63073}}
    \put(1010.12,229.751){\arc{90.0}{3.72512}{4.24872}}
    \put(1410.12,229.751){\arc{90.0}{3.72512}{4.24872}}
    \put(1369.88,310.249){\arc{90.0}{1.10713}{1.63073}}
    \put(1370.0,255.0){\arc{90.0}{4.18877}{4.71237}}
    \put(1370.0,345.0){\arc{90.0}{1.57078}{2.09438}}
    \put(1447.44,215.038){\arc{90.0}{3.20598}{3.72958}}
    \put(1372.56,264.962){\arc{90.0}{0.587983}{1.11158}}
    \put(1450.0,300.0){\circle*{12.0}}
    \put(1497.63,292.5){\arc{90.0}{3.66517}{4.18877}}
    \put(1450.0,350.0){\arc{50.0}{2.61797}{6.80674}}
    \path(1471.65,337.5)(1450.0,300.0)
    \path(1450.0,300.0)(1428.35,337.5)
    \put(1419.69,337.5){\arc{90.0}{0.0}{0.523579}}
    \path(1150.0,150.0)(1450.0,300.0)
    \path(1350.0,150.0)(1450.0,300.0)
    \path(1250.0,300.0)(1450.0,300.0)
    \path(1150.0,150.0)(1350.0,150.0)
    \path(1350.0,150.0)(1250.0,300.0)
    \path(1150.0,150.0)(1250.0,300.0)
    \put(1219.69,337.5){\arc{90.0}{0.0}{0.523579}}
    \path(1250.0,300.0)(1228.35,337.5)
    \path(1271.65,337.5)(1250.0,300.0)
    \put(1250.0,350.0){\arc{50.0}{2.61797}{6.80674}}
    \put(1297.63,292.5){\arc{90.0}{3.66517}{4.18877}}
    \put(1229.13,250.178){\arc{90.0}{5.69516}{6.21876}}
    \put(1304.02,300.101){\arc{90.0}{2.02997}{2.55357}}
    \put(1270.0,195.0){\arc{90.0}{1.57078}{2.09438}}
    \put(1270.0,105.0){\arc{90.0}{4.18877}{4.71237}}
    \put(1150.0,150.0){\circle*{12.0}}
    \put(1350.0,150.0){\circle*{12.0}}
    \put(1250.0,300.0){\circle*{12.0}}
    \put(1197.63,157.5){\arc{90.0}{2.09438}{2.61797}}
    \put(1150.0,100.0){\arc{50.0}{5.75957}{9.94834}}
    \path(1171.65,112.5)(1150.0,150.0)
    \path(1150.0,150.0)(1128.35,112.5)
    \put(1119.69,112.5){\arc{90.0}{5.75957}{6.28317}}
    \put(1397.63,157.5){\arc{90.0}{2.09438}{2.61797}}
    \put(1350.0,100.0){\arc{50.0}{5.75957}{9.94834}}
    \path(1371.65,112.5)(1350.0,150.0)
    \path(1350.0,150.0)(1328.35,112.5)
    \put(1319.69,112.5){\arc{90.0}{5.75957}{6.28317}}
    \put(1227.44,185.038){\arc{90.0}{3.72958}{4.25317}}
    \put(1152.56,234.962){\arc{90.0}{0.0643838}{0.587983}}
    \put(150.0,38.0){\makebox(0,0)[t]{$(a)$}}
    \put(500.0,38.0){\makebox(0,0)[t]{$(b)$}}
    \put(850.0,38.0){\makebox(0,0)[t]{$(c)$}}
    \put(1250.0,38.0){\makebox(0,0)[t]{$(d)$}}
\end{picture}
  \caption{Bases of $\alpha_2$, $\alpha_3$, $\alpha_4$ and $\zeta_4$}
  \label{nans.fig.small-twi}
\end{figure}
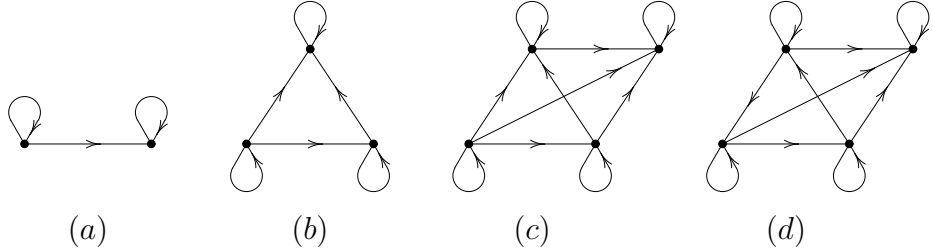

Let $D$ be a tournament with involution such that $|V(D)| \ge 2$. Let $S_1$, \ldots, $S_k$ be the
$\theta(D)$-classes of $D$. Recall that each $S_i$ takes the form $\{x, x'\}$ for some $x$.
Take arbitrary $x_1 \in S_1$. Then in each $S_j$, $j \ge 2$, one of the vertices dominates $x_1$ while
the other vertex is dominated by $x_1$. For each $j \in \{2, \ldots, k\}$ let $x_j \in S_j$ be the vertex
which dominates $x_1$. Clearly, $D[x_1, \ldots, x_k]$ is a reflexive torunament which, up to isomorphism,
uniquely determines~$D$. We shall say that $D[x_1, \ldots, x_k]$ is a \emph{base} of~$D$ and write
$D[x_1 \leftleftarrows x_2, \ldots, x_k]$ to emphasize the special status of $x_1$.
We say that a tournament with involution is \emph{acyclic} if each of its bases is an acyclic
reflexive tournament. Let $\alpha_n$ denote the acyclic tournament with involution with $2n$ vertices, $n \ge 1$,
and let $\alpha_0$ be the trivial one-vertex tournament with involution~$\1^\circ$.
The bases of $\alpha_2$, $\alpha_3$ and $\alpha_4$ are depicted in Fig.~\ref{nans.fig.small-twi} $(a)$, $(b)$ and
$(c)$, respectively. Let $\zeta_4$ denote the tournament with involution with $8$ vertices
whose base is depicted in Fig.~\ref{nans.fig.small-twi}~$(d)$.
Up to isomorphism, there are four distinct tournaments with involution with 4, 6 and 8 vertices:
$\alpha_2$, $\alpha_3$, $\alpha_4$ and $\zeta_4$,
and one can easily check that all of them are homomorphism-homogeneous.

\begin{LEM}\label{nans.lem.fuup}
  Let $D_1$ and $D_2$ be tournaments with involution and let $f : U \to W$ be a homomorphism from
  $D_1[U]$ to $D_2[W]$. Assume that there is a $u \in U$ such that $u' \in U$ and $f(u) = f(u')$.
  Then
  \begin{itemize}
  \item[$(a)$]
    $f(U) \subseteq \{v, v'\}$, where $v = f(u) = f(u')$;
  \item[$(b)$]
    $f$ extends to a homomorphism $f^*$ from $D_1$ to $D_2$.
  \end{itemize}
\end{LEM}
\begin{PROOF}
  $(a)$
  Take any $x \in U \setminus \{u, u'\}$. Since $D_1$ is a tournament with involution, we have
  $x \to u$ or $u \to x$, but not both. Say, $x \to u$. Then $u' \to x \to u$ whence
  $v \rightleftarrows f(x)$. Therefore, $f(x) \in \{v, v'\}$, since $D_2$ is also a tournament with involution.

  $(b)$ It is easy to see that $f^* : V(D_1) \to V(D_2)$ defined by $f^*(x) = x$ if $x \in U$ and
  $f^*(x) = v$ if $x \notin U$ is a homomorphism from $D_1$ to $D_2$.
\end{PROOF}

\begin{LEM}\label{nans.lem.fuup-2}
  Let $D_1$ and $D_2$ be tournaments with involution and let $f : U \to W$ be a homomorphism from
  $D_1[U]$ to $D_2[W]$. Assume that $f(U) \not \subseteq \{v, v'\}$ for all $v \in V(D_2)$.
  Then $f(u') = f(u)'$ whenever $u, u' \in U$.
\end{LEM}
\begin{PROOF}
  Assume that $u, u' \in U$. Since $u \rightleftarrows u'$, we have that $f(u) \rightleftarrows f(u')$,
  so $f(u') \in \{f(u), f(u)'\}$. If $f(u') = f(u)$ then, as we have just seen in Lemma~\ref{nans.lem.fuup},
  $f(U) \subseteq \{v, v'\}$ for $v = f(u)$, which is not the case. Therefore $f(u') = f(u)'$.
\end{PROOF}

\begin{LEM}\label{nans.lem.fuup-3}
  Let $D_1$ and $D_2$ be tournaments with involution and let $f : U \to W$ be a homomorphism from
  $D_1[U]$ to $D_2[W]$ such that $f(u') = f(u)'$ whenever $u, u' \in U$. Then $f$ extends to a
  homomorphism $\overline f$ from $D_1[\overline U]$ to $D_2[\overline W]$ where
  $\overline U = U \union \{u' : u \in U\}$ and $\overline W = W \union \{w' : w \in W\}$.
\end{LEM}
\begin{PROOF}
  Assume that there is an $x \in U$ such that $x' \notin U$. Then $f_1 : U \union \{x\} \to
  W \union \{f(x)'\}$ defined by
  $$
    f_1(u) = \begin{cases}
      f(u), & u \in U\\
      f(x)', & u = x'
    \end{cases}
  $$
  is a homomorphism from $D_1[U \union \{x\}]$ to $D_2[W \union \{f(x)'\}]$. We can repeat this
  procedure for every $x \in U$ such that $x' \notin U$ and thus extend $f$ to $\overline U$.
\end{PROOF}

\begin{LEM}\label{nans.lem.alpha-alpha}
  $\alpha_m \Rightarrow \alpha_n$ for all $m, n \ge 0$.
\end{LEM}
\begin{PROOF}
  If $m \le 1$ or $n \le 1$ the claim is trivially true. Assume, therefore, that $m, n \ge 2$.
  Fix a base $\alpha_m[x_1 \leftleftarrows x_2, \ldots, x_m]$ of $\alpha_m$.
  Let $f : U \to W$ be a
  homomorphism from $\alpha_m[U]$ to $\alpha_n[W]$ where $U = \{u_1, \ldots, u_l\}$ and
  $u_1 \to u_2 \to \ldots \to u_l$.
  Due to Lemmas~\ref{nans.lem.fuup} and \ref{nans.lem.fuup-2} without loss of generality we can assume that
  $f(u') = f(u)'$ whenever $u, u' \in U$.

  Let $\overline U = U \union \{x_1, \ldots, x_m\}$ and let $\overline f : \overline U \to W$
  be the mapping defined by $\overline f(u) = f(u)$ for $u \in U$ and
  $\overline f(x_i) = f(u_{\alpha(i)})$, $1 \le i \le m$, where
  $$
    \alpha(i) = \begin{cases}
      \min\{k : x_i \to u_k\}, & x_i \to u_l\\
      l, & x \not\to u_l.
    \end{cases}
  $$
  Then $\overline f$ is well-defined (if $x_i \in U$, say $x_i = u_j$ for some $j$, then
  $\overline f(x_i) = f(u_{\alpha(i)}) = f(u_j) = f(x_i)$) and it is a homomorphism from
  $\alpha_m[\overline U]$ to $\alpha_n[W]$ since $x_i \to x_j$ implies $u_{\alpha(i)} \to
  u_{\alpha(j)}$. According to Lemma~\ref{nans.lem.fuup-3}, $\overline f$ now easily extends
  to a homomorphism $f^* : \alpha_m \to \alpha_n$.
\end{PROOF}

\begin{figure}
  \centering
  \unitlength 0.24pt
\begin{picture}(1450.0,650.0)
    \path(178.054,122.746)(471.946,122.746)
    \put(325.0,575.0){\circle*{12.0}}
    \put(87.2359,402.254){\circle*{12.0}}
    \put(178.054,122.746){\circle*{12.0}}
    \put(471.946,122.746){\circle*{12.0}}
    \put(562.764,402.254){\circle*{12.0}}
    \put(372.631,567.5){\arc{90.0}{3.66517}{4.18877}}
    \put(325.0,625.0){\arc{50.0}{2.61797}{6.80674}}
    \path(346.651,612.5)(325.0,575.0)
    \path(325.0,575.0)(303.349,612.5)
    \put(294.689,612.5){\arc{90.0}{0.0}{0.523579}}
    \put(42.2047,385.015){\arc{90.0}{5.02653}{5.55013}}
    \path(87.2359,402.254)(44.8808,393.251)
    \path(58.2617,434.433)(87.2359,402.254)
    \put(39.683,417.705){\arc{50.0}{1.36134}{5.55013}}
    \put(109.088,445.237){\arc{90.0}{2.40853}{2.93213}}
    \put(180.534,74.5913){\arc{90.0}{3.76989}{4.29349}}
    \path(178.054,122.746)(173.527,79.6817)
    \path(138.496,105.134)(178.054,122.746)
    \put(148.664,82.2949){\arc{50.0}{0.1047}{4.29349}}
    \put(143.927,156.81){\arc{90.0}{1.1519}{1.6755}}
    \put(518.51,110.224){\arc{90.0}{2.51325}{3.03685}}
    \path(471.946,122.746)(511.504,105.134)
    \path(476.473,79.6817)(471.946,122.746)
    \put(501.336,82.2949){\arc{50.0}{5.13125}{9.32002}}
    \put(429.003,100.816){\arc{90.0}{6.17845}{6.70202}}
    \put(589.062,442.67){\arc{90.0}{1.25662}{1.78022}}
    \path(562.764,402.254)(591.738,434.433)
    \path(605.119,393.251)(562.764,402.254)
    \put(610.317,417.705){\arc{50.0}{3.87461}{8.06338}}
    \put(570.35,354.636){\arc{90.0}{4.92181}{5.44541}}
    \put(352.673,575.734){\makebox(0,0)[l]{$x_1$}}
    \path(562.764,402.254)(87.2359,402.254)
    \path(178.054,122.746)(325.0,575.0)
    \path(471.946,122.746)(325.0,575.0)
    \path(325.0,575.0)(87.2359,402.254)
    \path(325.0,575.0)(562.764,402.254)
    \put(494.109,507.758){\arc{90.0}{2.19909}{2.72269}}
    \put(441.208,434.947){\arc{90.0}{4.81709}{5.34069}}
    \put(208.792,434.947){\arc{90.0}{4.08405}{4.60765}}
    \put(155.891,507.758){\arc{90.0}{0.418859}{0.942458}}
    \put(252.813,498.455){\arc{90.0}{0.314139}{0.837738}}
    \put(338.408,470.643){\arc{90.0}{2.93213}{3.45573}}
    \put(311.592,470.643){\arc{90.0}{5.96901}{6.49258}}
    \put(397.187,498.455){\arc{90.0}{2.30381}{2.82741}}
    \path(178.054,122.746)(87.2359,402.254)
    \path(562.764,402.254)(471.946,122.746)
    \put(80.7655,276.545){\arc{90.0}{5.96901}{6.49258}}
    \put(166.361,304.357){\arc{90.0}{2.30381}{2.82741}}
    \put(354.389,167.746){\arc{90.0}{1.57078}{2.09438}}
    \put(354.389,77.7458){\arc{90.0}{4.18877}{4.71237}}
    \put(551.071,220.643){\arc{90.0}{3.45573}{3.97933}}
    \put(465.476,248.455){\arc{90.0}{6.07373}{6.5973}}
    \put(182.342,357.254){\arc{90.0}{4.71237}{5.23597}}
    \put(182.342,447.254){\arc{90.0}{1.04718}{1.57078}}
    \put(73.5741,378.953){\makebox(0,0)[tr]{$p'$}}
    \put(576.734,380.301){\makebox(0,0)[tl]{$s'$}}
    \put(198.882,101.175){\makebox(0,0)[tl]{$q$}}
    \put(452.777,102.523){\makebox(0,0)[tr]{$r$}}
    \put(1252.78,102.523){\makebox(0,0)[tr]{$r$}}
    \put(998.882,101.175){\makebox(0,0)[tl]{$q$}}
    \put(1376.73,380.301){\makebox(0,0)[tl]{$s$}}
    \put(873.574,378.953){\makebox(0,0)[tr]{$p$}}
    \put(1154.39,77.7458){\arc{90.0}{4.18877}{4.71237}}
    \put(1154.39,167.746){\arc{90.0}{1.57078}{2.09438}}
    \put(1197.19,498.455){\arc{90.0}{2.30381}{2.82741}}
    \put(1111.59,470.643){\arc{90.0}{5.96901}{6.49258}}
    \put(1138.41,470.643){\arc{90.0}{2.93213}{3.45573}}
    \put(1052.81,498.455){\arc{90.0}{0.314139}{0.837738}}
    \path(1271.95,122.746)(1125.0,575.0)
    \path(978.054,122.746)(1125.0,575.0)
    \put(1152.67,575.734){\makebox(0,0)[l]{$x_1$}}
    \put(1370.35,354.636){\arc{90.0}{4.92181}{5.44541}}
    \put(1410.32,417.705){\arc{50.0}{3.87461}{8.06338}}
    \path(1405.12,393.251)(1362.76,402.254)
    \path(1362.76,402.254)(1391.74,434.433)
    \put(1389.06,442.67){\arc{90.0}{1.25662}{1.78022}}
    \put(1229.0,100.816){\arc{90.0}{6.17845}{6.70202}}
    \put(1301.34,82.2949){\arc{50.0}{5.13125}{9.32002}}
    \path(1276.47,79.6817)(1271.95,122.746)
    \path(1271.95,122.746)(1311.5,105.134)
    \put(1318.51,110.224){\arc{90.0}{2.51325}{3.03685}}
    \put(943.927,156.81){\arc{90.0}{1.1519}{1.6755}}
    \put(948.664,82.2949){\arc{50.0}{0.1047}{4.29349}}
    \path(938.496,105.134)(978.054,122.746)
    \path(978.054,122.746)(973.527,79.6817)
    \put(980.534,74.5913){\arc{90.0}{3.76989}{4.29349}}
    \put(909.088,445.237){\arc{90.0}{2.40853}{2.93213}}
    \put(839.683,417.705){\arc{50.0}{1.36134}{5.55013}}
    \path(858.262,434.433)(887.236,402.254)
    \path(887.236,402.254)(844.881,393.251)
    \put(842.205,385.015){\arc{90.0}{5.02653}{5.55013}}
    \put(1094.69,612.5){\arc{90.0}{0.0}{0.523579}}
    \path(1125.0,575.0)(1103.35,612.5)
    \path(1146.65,612.5)(1125.0,575.0)
    \put(1125.0,625.0){\arc{50.0}{2.61797}{6.80674}}
    \put(1172.63,567.5){\arc{90.0}{3.66517}{4.18877}}
    \put(1362.76,402.254){\circle*{12.0}}
    \put(1271.95,122.746){\circle*{12.0}}
    \put(978.054,122.746){\circle*{12.0}}
    \put(887.236,402.254){\circle*{12.0}}
    \put(1125.0,575.0){\circle*{12.0}}
    \path(978.054,122.746)(1271.95,122.746)
    \path(887.236,402.254)(978.054,122.746)
    \path(1271.95,122.746)(887.236,402.254)
    \path(887.236,402.254)(1362.76,402.254)
    \path(978.054,122.746)(1362.76,402.254)
    \path(1271.95,122.746)(1362.76,402.254)
    \put(984.524,248.455){\arc{90.0}{2.82741}{3.35101}}
    \put(898.929,220.643){\arc{90.0}{5.44541}{5.96901}}
    \put(937.728,309.947){\arc{90.0}{5.34069}{5.86429}}
    \put(990.628,382.758){\arc{90.0}{1.6755}{2.19909}}
    \put(1125.0,447.254){\arc{90.0}{1.57078}{2.09438}}
    \put(1125.0,357.254){\arc{90.0}{4.18877}{4.71237}}
    \put(1143.96,298.906){\arc{90.0}{0.942458}{1.46606}}
    \put(1196.86,226.094){\arc{90.0}{3.56045}{4.08405}}
    \put(1274.56,276.406){\arc{90.0}{0.314139}{0.837738}}
    \put(1360.15,248.594){\arc{90.0}{2.93213}{3.45573}}
    \path(887.236,402.254)(1125.0,575.0)
    \path(1362.76,402.254)(1125.0,575.0)
    \put(1193.66,469.496){\arc{90.0}{5.34069}{5.86429}}
    \put(1246.56,542.307){\arc{90.0}{1.6755}{2.19909}}
    \put(1003.44,542.307){\arc{90.0}{0.942458}{1.46606}}
    \put(1056.34,469.496){\arc{90.0}{3.56045}{4.08405}}
    \put(325.0,38.0){\makebox(0,0)[t]{$(a)$}}
    \put(1125.0,38.0){\makebox(0,0)[t]{$(b)$}}
\end{picture}
  \caption{The proof of Theorem~\ref{nans.thm.hh-tinv}}
  \label{nans.fig.acyc-c4}
\end{figure}
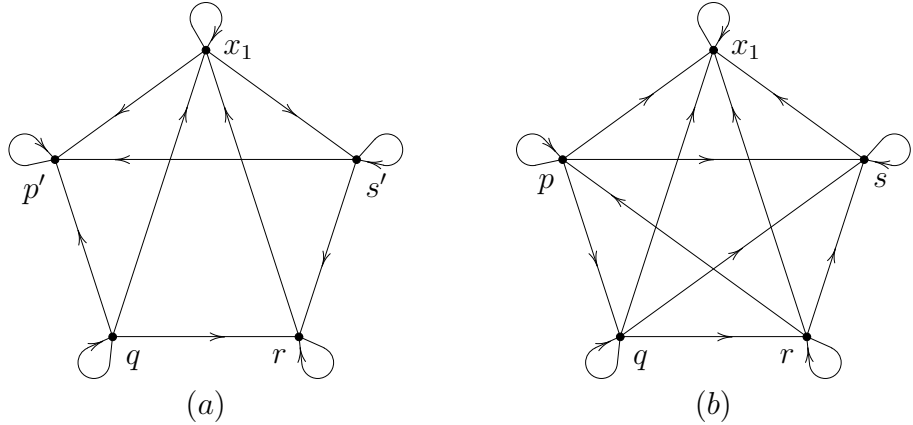

\begin{THM}\label{nans.thm.hh-tinv}
  Let $D$ be a tournament with involution. Then $D$ is homo\-morph\-ism-homogeneous if and only if
  \begin{enumerate}
  \item
    $D \cong \zeta_4$, or
  \item
    $D \cong \alpha_n$ for some $n \ge 0$.
  \end{enumerate}
\end{THM}
\begin{PROOF}
  $(\Leftarrow)$
  We have already seen that $\zeta_4$ is homomorphism-homogeneous.
  From Lemma~\ref{nans.lem.alpha-alpha} it follows that
  $\alpha_n \Rightarrow \alpha_n$ for all $n \ge 1$, so $\alpha_n$ is homo\-morph\-ism-homogeneous for all~$n \ge 0$.

  $(\Rightarrow)$
  Let $D$ be a homomorphism-homogeneous tournament with involution.
  If $|V(D)| \le 8$ then $D \cong \zeta_4$ or $D \cong \alpha_n$ for some $n \in \{0, 1, 2, 3, 4\}$.
  Assume now that $|V(D)| \ge 10$ and let us show that
  $D$ is an acyclic tournament with involution.

  Suppose, to the contrary, that $D$ is not an acyclic tournament with involution
  and let $D[x_1 \leftleftarrows x_2, \ldots, x_k]$
  be a base od $D$ which is not an acyclic tournament. Then $D[x_2, \ldots, x_k]$ is not acyclic.
  Since $|V(D)| \ge 10$, every base of $D$ has at least 5 vertices, so $k \ge 5$.

  Assume that there exist distinct $y_1, \ldots, y_m \in \{x_2, \ldots, x_k\}$ such that
  $y_1 \to y_2 \to \ldots \to y_m \to y_1$ is a cycle of length $m \ge 4$. Then $D[y_1, \ldots, y_m]$
  is a Hamiltonian tournament and hence pancyclic. Therefore, there exist four distinct vertices
  $p, q, r, s \in \{y_1, \ldots, y_m\}$ such that $p \to q \to r \to s \to p$ is a 4-cycle,
  Fig.~\ref{nans.fig.acyc-c4}~$(a)$ (instead of $p$ and $q$ the figure depicts $p'$ and $q'$).
  The mapping
  $$
    f = \begin{pmatrix}
      p' & q & r & s'\\
      r  & q & r & q
    \end{pmatrix}
  $$
  is a homomorphism from $D[p', q, r, s']$ to $D[q, r]$, so it extends to an endomorphism $f^*$ of $D$.
  From $r \to x_1 \to p'$ we have $f^*(r) \to f^*(x_1) \to f^*(p')$, that is, $r \to f^*(x_1) \to r$.
  Therefore, $f^*(x_1) \in \{r, r'\}$. Analogously, from $q \to x_1 \to s'$ we conclude
  $f^*(x_1) \in \{q, q'\}$. Contradiction.

  Assume now that there exist distinct $p, q, r \in \{x_2, \ldots, x_k\}$ such that
  $p \to q \to r \to p$ is a 3-cycle. Since $k \ge 5$ there exists an $s \in \{x_2, \ldots, x_k\}
  \setminus \{p, q, r\}$. As we have just seen, $D[x_2, \ldots, x_k]$ does not contain a 4-cycle,
  so $\{p, q, r\} \TO s$ or $s \TO \{p, q, r\}$. Without loss of generality we can assume that
  $\{p, q, r\} \TO s$, Fig.~\ref{nans.fig.acyc-c4}~$(b)$.
  The mapping
  $$
    f = \begin{pmatrix}
      p & q & s & x_1\\
      p & q & q & x_1
    \end{pmatrix}
  $$
  is a homomorphism from $D[p, q, s, x_1]$ to $D[p, q, x_1]$, so it extends to an endomorphism $f^*$ of $D$.
  From $q \to r \to s$ we conclude that $q \to f^*(r) \to q$, so $f^*(r) \in \{q, q'\}$. If $f^*(r) = q$ then
  $r \to p$ implies $q \to p$, which is not the case. On the other hand, if $f^*(r) = q'$ then
  $r \to x_1$ implies $q' \to x_1$, which is not possible (since $q \to x_1$ enforces $x_1 \to q'$).
  Therefore, $D[x_2, \ldots, x_k]$ is acyclic, and hence $D[x_1 \leftleftarrows x_2, \ldots, x_k]$
  is acyclic.
\end{PROOF}

\begin{figure}
  \centering
  \unitlength 0.24pt
\begin{picture}(1300.0,375.0)
    \path(500.0,150.0)(700.0,150.0)
    \put(469.689,187.5){\arc{90.0}{0.0}{0.523579}}
    \path(500.0,150.0)(478.349,187.5)
    \path(521.651,187.5)(500.0,150.0)
    \put(500.0,200.0){\arc{50.0}{2.61797}{6.80674}}
    \put(547.631,142.5){\arc{90.0}{3.66517}{4.18877}}
    \put(669.689,187.5){\arc{90.0}{0.0}{0.523579}}
    \path(700.0,150.0)(678.349,187.5)
    \path(721.651,187.5)(700.0,150.0)
    \put(700.0,200.0){\arc{50.0}{2.61797}{6.80674}}
    \put(747.631,142.5){\arc{90.0}{3.66517}{4.18877}}
    \put(620.0,195.0){\arc{90.0}{1.57078}{2.09438}}
    \put(620.0,105.0){\arc{90.0}{4.18877}{4.71237}}
    \put(500.0,150.0){\circle*{12.0}}
    \put(700.0,150.0){\circle*{12.0}}
    \put(310.125,229.751){\arc{90.0}{3.72512}{4.24872}}
    \put(269.875,310.249){\arc{90.0}{1.10713}{1.63073}}
    \put(270.0,255.0){\arc{90.0}{4.18877}{4.71237}}
    \put(270.0,345.0){\arc{90.0}{1.57078}{2.09438}}
    \put(347.442,215.038){\arc{90.0}{3.20598}{3.72958}}
    \put(272.558,264.962){\arc{90.0}{0.587983}{1.11158}}
    \put(350.0,300.0){\circle*{12.0}}
    \put(397.631,292.5){\arc{90.0}{3.66517}{4.18877}}
    \put(350.0,350.0){\arc{50.0}{2.61797}{6.80674}}
    \path(371.651,337.5)(350.0,300.0)
    \path(350.0,300.0)(328.349,337.5)
    \put(319.689,337.5){\arc{90.0}{0.0}{0.523579}}
    \path(50.0,150.0)(350.0,300.0)
    \path(250.0,150.0)(350.0,300.0)
    \path(150.0,300.0)(350.0,300.0)
    \path(50.0,150.0)(250.0,150.0)
    \path(250.0,150.0)(150.0,300.0)
    \path(50.0,150.0)(150.0,300.0)
    \put(119.689,337.5){\arc{90.0}{0.0}{0.523579}}
    \path(150.0,300.0)(128.349,337.5)
    \path(171.651,337.5)(150.0,300.0)
    \put(150.0,350.0){\arc{50.0}{2.61797}{6.80674}}
    \put(197.631,292.5){\arc{90.0}{3.66517}{4.18877}}
    \put(129.132,250.178){\arc{90.0}{5.69516}{6.21876}}
    \put(204.016,300.101){\arc{90.0}{2.02997}{2.55357}}
    \put(170.0,195.0){\arc{90.0}{1.57078}{2.09438}}
    \put(170.0,105.0){\arc{90.0}{4.18877}{4.71237}}
    \put(50.0,150.0){\circle*{12.0}}
    \put(250.0,150.0){\circle*{12.0}}
    \put(150.0,300.0){\circle*{12.0}}
    \put(97.6314,157.5){\arc{90.0}{2.09438}{2.61797}}
    \put(50.0,100.0){\arc{50.0}{5.75957}{9.94834}}
    \path(71.6506,112.5)(50.0,150.0)
    \path(50.0,150.0)(28.3494,112.5)
    \put(19.6891,112.5){\arc{90.0}{5.75957}{6.28317}}
    \put(297.631,157.5){\arc{90.0}{2.09438}{2.61797}}
    \put(250.0,100.0){\arc{50.0}{5.75957}{9.94834}}
    \path(271.651,112.5)(250.0,150.0)
    \path(250.0,150.0)(228.349,112.5)
    \put(219.689,112.5){\arc{90.0}{5.75957}{6.28317}}
    \put(127.442,185.038){\arc{90.0}{3.72958}{4.25317}}
    \put(52.5577,234.962){\arc{90.0}{0.0643838}{0.587983}}
    \put(500.0,132.0){\makebox(0,0)[t]{$u$}}
    \put(700.0,132.0){\makebox(0,0)[t]{$t$}}
    \put(37.9783,162.453){\makebox(0,0)[br]{$p$}}
    \put(268.0,150.0){\makebox(0,0)[l]{$q$}}
    \put(132.0,300.0){\makebox(0,0)[r]{$r$}}
    \put(362.384,287.904){\makebox(0,0)[tl]{$s$}}
    \put(600.0,38.0){\makebox(0,0)[t]{$\alpha_2$}}
    \put(150.0,38.0){\makebox(0,0)[t]{$\zeta_4$}}
    \put(1210.12,229.751){\arc{90.0}{3.72512}{4.24872}}
    \put(1169.88,310.249){\arc{90.0}{1.10713}{1.63073}}
    \put(1170.0,255.0){\arc{90.0}{4.18877}{4.71237}}
    \put(1170.0,345.0){\arc{90.0}{1.57078}{2.09438}}
    \put(1247.44,215.038){\arc{90.0}{3.20598}{3.72958}}
    \put(1172.56,264.962){\arc{90.0}{0.587983}{1.11158}}
    \put(1250.0,300.0){\circle*{12.0}}
    \put(1297.63,292.5){\arc{90.0}{3.66517}{4.18877}}
    \put(1250.0,350.0){\arc{50.0}{2.61797}{6.80674}}
    \path(1271.65,337.5)(1250.0,300.0)
    \path(1250.0,300.0)(1228.35,337.5)
    \put(1219.69,337.5){\arc{90.0}{0.0}{0.523579}}
    \path(950.0,150.0)(1250.0,300.0)
    \path(1150.0,150.0)(1250.0,300.0)
    \path(1050.0,300.0)(1250.0,300.0)
    \path(950.0,150.0)(1050.0,300.0)
    \put(1019.69,337.5){\arc{90.0}{0.0}{0.523579}}
    \path(1050.0,300.0)(1028.35,337.5)
    \path(1071.65,337.5)(1050.0,300.0)
    \put(1050.0,350.0){\arc{50.0}{2.61797}{6.80674}}
    \put(1097.63,292.5){\arc{90.0}{3.66517}{4.18877}}
    \put(972.558,264.962){\arc{90.0}{0.587983}{1.11158}}
    \put(1047.44,215.038){\arc{90.0}{3.20598}{3.72958}}
    \put(950.0,150.0){\circle*{12.0}}
    \put(1150.0,150.0){\circle*{12.0}}
    \put(1050.0,300.0){\circle*{12.0}}
    \put(997.631,157.5){\arc{90.0}{2.09438}{2.61797}}
    \put(950.0,100.0){\arc{50.0}{5.75957}{9.94834}}
    \path(971.651,112.5)(950.0,150.0)
    \path(950.0,150.0)(928.349,112.5)
    \put(919.689,112.5){\arc{90.0}{5.75957}{6.28317}}
    \put(1197.63,157.5){\arc{90.0}{2.09438}{2.61797}}
    \put(1150.0,100.0){\arc{50.0}{5.75957}{9.94834}}
    \path(1171.65,112.5)(1150.0,150.0)
    \path(1150.0,150.0)(1128.35,112.5)
    \put(1119.69,112.5){\arc{90.0}{5.75957}{6.28317}}
    \put(937.978,162.453){\makebox(0,0)[br]{$u$}}
    \put(1168.0,150.0){\makebox(0,0)[l]{$x_i$}}
    \put(1262.38,287.904){\makebox(0,0)[tl]{$t$}}
    \put(1032.0,300.0){\makebox(0,0)[r]{$v$}}
    \put(1050.0,38.0){\makebox(0,0)[t]{$\alpha_n$}}
    \dottedline{12}(1050.0,300.0)(1150.0,150.0)
    \dottedline{12}(1150.0,150.0)(950.0,150.0)
    \put(1010.0,105.0){\arc{90.0}{4.71237}{5.23597}}
    \put(1010.0,195.0){\arc{90.0}{1.04718}{1.57078}}
    \put(1157.44,219.962){\arc{90.0}{2.55357}{3.07717}}
    \put(1082.56,170.038){\arc{90.0}{5.17156}{5.69516}}
\end{picture}
  \caption{The proof of Lemma~\ref{nans.lem.zeta-alpha}}
  \label{nans.fig.zeta-alpha}
\end{figure}
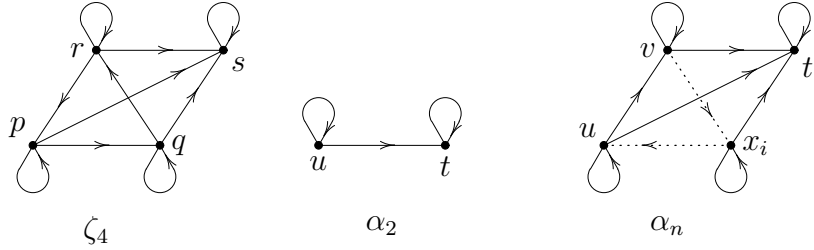

\begin{LEM}\label{nans.lem.zeta-alpha}
  For all $n \ge 2$ we have $\zeta_4 \not\Rightarrow \alpha_n$.
\end{LEM}
\begin{PROOF}
  Assume first that $\zeta_4 \Rightarrow \alpha_2$.
  Let $\zeta_4[s \leftleftarrows p, q, r]$ be a base of $\zeta_4$ such that $p \to q \to r \to p$, and
  let $\alpha_2[t \leftleftarrows u]$ be a base of $\alpha_2$, Fig.~\ref{nans.fig.zeta-alpha}. The mapping
  $$
    f = \begin{pmatrix}
      p & s & q \\
      u & t & t
    \end{pmatrix}
  $$
  is a homomorphism from $\zeta_4[s, p, q]$ to $\alpha_2[t, u]$, so by the assumption it extends to
  a homomorphism $f^*$ from $\zeta_4$ to $\alpha_2$. Let us compute $f^*(r)$. From $r \to s$ it follows
  that $f^*(r) \to t$. Therefore, $f^*(r)$ is a vertex in the base of $\alpha_2$ under consideration,
  so $f^*(r) \in \{u, t\}$. If $f^*(r) = u$ then $q \to r$ implies $f^*(q) = t \to u = f^*(r)$, which is not the
  case. On the other hand, if $f^*(r) = t$ then $r \to p$ implies $f^*(r) = t \to u = f^*(p)$, which is not the
  case. This contradiction shows that $\zeta_4 \not\Rightarrow \alpha_2$.

  Assume now that $\zeta_4 \Rightarrow \alpha_n$ for some $n \ge 3$.
  As above, let $\zeta_4[s \leftleftarrows p, q, r]$ be a base of $\zeta_4$ such that $p \to q \to r \to p$, and
  let $\alpha_n[t \leftleftarrows u, v, x_4, \ldots, x_n]$ be a base of $\alpha_n$ such that $u \to v$.
  The mapping
  $$
    f = \begin{pmatrix}
      p & s & q \\
      u & t & v
    \end{pmatrix}
  $$
  is a homomorphism from $\zeta_4[s, p, q]$ to $\alpha_n[t, u, v]$, so by the assumption it extends to
  a homomorphism $f^*$ from $\zeta_4$ to $\alpha_n$. Let us compute $f^*(r)$. From $r \to s$ it follows
  that $f^*(r) \to t$. Therefore, $f^*(r)$ is a vertex in the base of $\alpha_n$ under consideration,
  that is, $f^*(r) \in \{t, u, v, x_4, \ldots, x_n\}$. If $f^*(r) = t$ then $r \to p$ implies
  $f^*(r) = t \to u = f^*(p)$, which is not the case.  If $f^*(r) = u$ then $q \to r$ implies
  $f^*(q) = v \to u = f^*(p)$, which is not the case.  If $f^*(r) = v$ then $r \to p$ implies
  $f^*(r) = v \to u = f^*(p)$, which is not the case.  Finally, if $f^*(r) = x_i$ for some $i$
  then $q \to r \to p$ implies $v \to x_i \to u$, which is not the case since
  $\alpha_n[t \leftleftarrows u, v, x_4, \ldots, x_n]$ is an acyclic digraph.
  Therefore, $\zeta_4 \not\Rightarrow \alpha_n$.
\end{PROOF}

\begin{THM}\label{nans.thm.hh-dwi}
  Let $D$ be a digraph with involution. Then $D$ is homomorphism-homogeneous if and only if
  \begin{enumerate}
  \item
    $D \cong k \cdot \alpha_0 + l \cdot \alpha_1 + m \cdot \zeta_4$ for some $k, l, m \ge 0$, or
  \item
    $D \cong m_1 \cdot \alpha_{n_1} + \ldots + m_k \cdot \alpha_{n_k}$
    for some $k \ge 1$ and $m_1, \ldots, m_k, n_1, \ldots, n_k \ge 0$.
  \end{enumerate}
\end{THM}
\begin{PROOF}
  $(\Leftarrow)$
  It is easy to see that $D_1 \Rightarrow D_2$ for all choices of $D_1, D_2 \in \{\alpha_0, \alpha_1, \zeta_4\}$,
  so digraphs from the class (1) are homomorphism-homogeneous. We have seen in Lemma~\ref{nans.lem.alpha-alpha}
  that $\alpha_{n_i} \Rightarrow \alpha_{n_j}$ for all $n_i, n_j \ge 0$, so digraphs
  from the class (2) are also homomorphism-homogeneous.

  $(\Rightarrow)$
  Let $D$ be a homomorphism-homogeneous digraph with involution.
  Then by Lemma~\ref{nans.lem.dprime-hh} every connected component of $D$ is a homomorphism-homogeneous
  tournament with involution. Therefore, Theorem~\ref{nans.thm.hh-tinv} yields that every connected
  component of $D$ is isomorphic to $\zeta_4$ or $\alpha_n$ for some $n \ge 0$.
  If there is a connected component of $D$ isomorphic to $\zeta_4$, then due to Lemma~\ref{nans.lem.zeta-alpha}
  every connected component of $D$ is isomorphic to $\zeta_4$, $\alpha_0$ or $\alpha_1$, and we have
  case (1). On the other hand, if no connected component of $D$ is isomorphic to $\zeta_4$ then
  every connected component of $D$ is isomorphic to $\alpha_n$ for some $n \ge 0$ and we have case (2).
\end{PROOF}

\section{Bidirectionally disconnected systems}
\label{nans.sec.bi-disconn}

Reflexive homomorphism-homogeneous proper digraphs were characterized in \cite[Theorem 3.10]{masul-hhd},
see Theorem~\ref{nans.thm.hh-digraph}. As for bidirectionally connected
digraphs, we have seen in Theorem~\ref{nans.thm.KK1} that
we cannot hope for a reasonable description due to the complexity of the
corresponding decision problem.
In this section we characterize finite reflexive homomorphism-homogeneous bidirectionally disconnected
improper digraphs.

Let us start with a rather general result. We say that a digraph $(V^*, E^*)$ is a
\emph{retract} of a digraph $(V, E)$ if there exist homomorphisms
$r : V \to V^*$ and $j : V^* \to V$ such that $r \circ j = \id_{V^*}$.

\begin{LEM}\label{nans.lem.al-retr}
  Let $D$ be a reflexive improper digraph and let $\rho \subseteq V(D)^2$ be an equivalence
  relation on $V(D)$ such that the following holds:
  \begin{itemize}
  \item
    for every $S \in V(D) / \rho$ we have $D[S] \cong K_n^\circ$ for some positive integer $n$;
  \item
    for all distinct $S$, $T \in V(D) / \rho$, if $S \sim T$ then $S \TO T$ or
    $T \TO S$ or both.
  \end{itemize}
  Define the digraph $D / \rho$ as follows: the set of vertices of $D / \rho$ is
  $V(D) / \rho$, while $(S, T)$ is an edge od $D / \rho$ if and only if $S = T$ or
  $S \TO T$ in $D$. Then
  \begin{enumerate}
  \item
    $D / \rho$ is a reflexive digraph.
  \item
    $D / \rho$ is a retract of $D$.
  \item
    If $D$ is a homomorphism-homogeneous digraph then $D / \rho$ is
    a homomorphism-homogeneous digraph.
  \item
    Assume that the following holds for $D$:
    \begin{itemize}
    \item[$(\diamondsuit)$]
      for all distinct $S$, $T \in V(D) / \rho$, if $S \sim T$ then $S \TO T$ or
      $T \TO S$, but not both.
    \end{itemize}
    Then $D$ is a homomorphism-homogeneous digraph if and only if $D / \rho$ is
    a homomorphism-homogeneous digraph.
  \end{enumerate}
\end{LEM}
\begin{PROOF}
  Let $V(D) / \rho = \{S_1, \ldots, S_n\}$.

  (1) Obvious.

  (2) Choose arbitrary $s_1 \in S_1$, \ldots, $s_n \in S_n$ and define $r : V(D) \to V(D/\rho)$ and
  $j : V(D/\rho) \to V(D)$ by $r(x) = x / \rho$ and $j(S_i) = s_i$.
  Then $r$ and $j$ are homomorphisms satisfying $r \circ j = \id_{V(D/\rho)}$, so $D / \rho$
  is a retract of~$D$.

  (3) It is easy to see that every retract of a homomorphism-homogeneous relational
  structure is homomorphism-homogeneous.
  Therefore, $D / \rho$, being a retract of $D$, is homomorphism-homogeneous.

  (4) Direction from left to right follows from (3). Let us show the other direction.
  Let $f : U \to W$ be a homomorphism from $D[U]$ to $D[W]$ where $U, W \subseteq V(D)$.
  From $(\diamondsuit)$ it follows that if $U \sec S \ne \0$ for some $S \in V(D)/\rho$ then
  $f(U \sec S) \subseteq S'$ for some $S' \in V(D)/\rho$. Without loss of generality,
  let $S_1, \ldots, S_k$ be all the
  $\rho$-classes that intersect $U$ and let $S'_1, \ldots, S'_k$ be the $\rho$-classes such that
  $f(U \sec S_i) \subseteq S'_i$, $1 \le i \le k$. Then the mapping
  $g : \{S_1, \ldots, S_k\} \to \{S'_1, \ldots, S'_k\} : S_i \mapsto S'_i$ is easily seen to be a
  homomorphism from $(D/\rho)[S_1, \ldots, S_k]$ to $(D/\rho)[S'_1, \ldots, S'_k]$.
  Since $D/\rho$ is homomorphism-homogeneous, $g$ extends to an endomorphism $g^*$ of $D/\rho$.
  Then the mapping $f^* : V(D) \to V(D)$ defined by
  $$
    f^*(x) = \begin{cases}
      f(x), & x \in U\\
      j \circ g^* \circ r(x), & x \notin U
    \end{cases}
  $$
  where $r$ and $j$ are the homomorphisms from (2), is an endomorphism of $D$ which extends~$f$.
\end{PROOF}

Recall that a quasiorder is a binary relational system $(A, \le)$ where $\le$ is a reflexive and transitive binary
relation on $A$. If $\equiv$ denotes the equivalence relation on $A$ defined by $x \equiv y$ if $x \le y$
and $y \le x$, then $A / \hbox{$\equiv$}$ is a partially ordered set where $x / \hbox{$\equiv$} \le
y / \hbox{$\equiv$}$ if and only if $x \le y$. Since $(A / \hbox{$\equiv$}, \le)$ is a retract of
$(A, \le)$, as a direct consequence of the above lemma we have the following:

\begin{COR}\label{nans.cor.quasiorder}
  Let $(A, \le)$ be a quasiorder. Then
  $(A, \le)$ is homomorphism-homogeneous as a quasiorder if and only if
  $(A / \hbox{$\equiv$}, \le)$ is a homomorphism-homogeneous partially ordered set.
\end{COR}

Let $D = (V, E)$ be a proper digraph with $V = \{ v_1, \ldots, v_n \}$, and let $V_1$, \ldots, $V_n$ be
finite nonempty pairwise disjoint sets. Let $D\langle V_1, \ldots, V_n \rangle$ denote the digraph
whose vertices  are $V_1 \union \ldots \union V_n$ and whose edges are defined as follows:
\begin{itemize}
\item
  for every $i \in \{1, \ldots, n\}$ and for all $x, y \in V_i$ we have $x \to y$ in $D\langle V_1, \ldots, V_n \rangle$;
\item
  if $v_i \to v_j$ in $D$ and $i \ne j$, then $V_i \TO V_j$ in $D\langle V_1, \ldots, V_n \rangle$;
\item
  no other edges exist in $D\langle V_1, \ldots, V_n \rangle$.
\end{itemize}
We say that $D\langle V_1, \ldots, V_n \rangle$ in an \emph{inflation} of $D$. Note that
$D\langle V_1, \ldots, V_n \rangle[V_i] \cong K^\circ_{|V_i|}$ and that $D$ is a retract of
$D\langle V_1, \ldots, V_n \rangle$.

\begin{LEM}\label{nans.lem.inflation}
  A finite reflexive proper digraph $D$ is homomorphism-homogeneous if and only if every inflation of $D$
  is homomorphism-homogeneous.
\end{LEM}

\begin{LEM}\label{nans.lem.al-1}
  Let $D$ be a finite homomorphism-homogeneous reflexive bidirectionally disconnected improper digraph,
  and let $S \in V(D) / \theta(D)$ be an arbitrary equivalence class of $\theta(D)$.
  Then $D[S] \cong K_n^\circ$ for some positive integer~$n$.
\end{LEM}
\begin{PROOF}
  From $\omega(D) < |V(D) / \theta(D)|$ it follows that there exist distinct classes of
  $\theta(D)$ which belong to the same connected component of $D$. Therefore, we can choose
  $T_1, T_2 \in V(D) / \theta(D)$ in such a way that $T_1 \ne T_2$ and $T_1 \to T_2$.
  Moreover, choose an $y_1 \in T_1$ and an $y_2 \in T_2$ so that $y_1 \to y_2$.

  Assume that there is an $S \in V(D) / \theta(D)$ such that $D[S]$ is not a complete reflexive graph.
  Then there exist $u, v \in S$ such that $u \not\rightleftarrows v$.
  If $u \not\sim v$ or $u \to v$, consider the mapping
  $$
    f : \begin{pmatrix}
      u   & v   \\
      y_1 & y_2
    \end{pmatrix}.
  $$
  If $v \to u$, consider
  $$
    f : \begin{pmatrix}
      v   & u   \\
      y_1 & y_2
    \end{pmatrix}.
  $$
  In any case, the mapping $f$ is a homomorphism between finitely induced subdigraphs of $D$, so it
  extends to an endomorphism $f^*$ of~$D$. Since $u$ and $v$ belong to
  the same equivalence class of $\theta(D)$, there exist $z_1, z_2, \ldots, z_k \in V(D)$ such that
  $$
    u = z_1 \rightleftarrows z_2 \rightleftarrows \ldots \rightleftarrows z_k = v.
  $$
  Then
  $
    f^*(u) = f^*(z_1) \rightleftarrows f^*(z_2) \rightleftarrows \ldots \rightleftarrows f^*(z_k) = f^*(v),
  $
  whence follows that $(y_1, y_2) \in \theta(D)$ since $\{f^*(u), f^*(v)\} = \{y_1, y_2\}$.
  Contradiction.
\end{PROOF}

Bidirectionally disconnected digraphs naturaly split into two classes:
\begin{itemize}
\item
  we say that a digraph $D$ is a digraph \emph{with no back-and-forth} if the following holds
  for all $S$, $T \in V(D) / \theta(D)$: if $S \rightleftarrows T$ then $S = T$;
\item
  we say that a digraph $D$ is a digraph \emph{with back-and-forth} if
  there exist distinct $S$, $T \in V(D) / \theta(D)$ such that $S \rightleftarrows T$.
\end{itemize}
Let us first classify homomorphism-homogeneous bidirectionally disconnected digraphs with
no back-and-forth.

\begin{LEM}\label{nans.lem.al-2}
  Let $D$ be a finite homomorphism-homogeneous reflexive bidirectionally disconnected improper digraph
  with no back-and-forth.
  Then for all distinct $S, T \in V(D) / \theta(D)$ either $S \not\sim T$, or $S \TO T$
  or $T \TO S$.
\end{LEM}
\begin{PROOF}
  Take any $S, T \in V(D) / \theta(D)$ such that $S \sim T$, assume that $S \to T$ and let
  us show that $S \TO T$. Assume, to the contrary, that there exist $v \in S$ and $w \in T$ such
  that $v \not\to w$.
  Since $S \ne T$ and $S \to T$, we know that $T \not\to S$, so $w \not\to v$ and thus
  $v \not\sim w$. Then the mapping
  $$
    f : \begin{pmatrix}
      v & w \\
      w & v
    \end{pmatrix}
  $$
  is a homomorphism between finitely induced subdigraphs of $D$, so it
  extends to an endomorphism $f^*$ of~$D$. Choose $x \in S$ and $y \in T$ so that $x \to y$.
  From $v \to x \to y \to w$ it follows that $f^*(v) \to f^*(x) \to f^*(y) \to f^*(w)$, that is,
  $w \to f^*(x) \to f^*(y) \to v$. Clearly, at least one of the edges
  $w \to f^*(x)$, $f^*(x) \to f^*(y)$ or $f^*(y) \to v$ leads from $T$ to $S$, which contradicts
  the fact that $T \not\to S$.
\end{PROOF}

\begin{THM}\label{nans.thm.case1a}
  Let $D$ be a finite reflexive bidirectionally disconnected improper digraph with no back-and-forth.
  Then $D$ is homomorphism-homo\-gene\-ous if and only if
  \begin{enumerate}
  \item
    $D$ is a finite homomorphism-homogeneous quasiorder; or
  \item
    $D$ is an inflation of
    $k \cdot C_3^\circ + l \cdot \1^\circ$ for some $k, l \ge 0$ such that
    $k + l \ge 1$.
  \end{enumerate}
\end{THM}
\begin{PROOF}
  Let $D$ be a finite reflexive bidirectionally disconnected improper digraph with no back-and-forth.

  $(\Rightarrow)$
  Assume that $D$ is a homomorphism-homogeneous digraph.
  According to Lemma~\ref{nans.lem.al-retr},  $D / \theta(D)$ is a homomorphism-homogeneous
  reflexive proper digraph, so Theorem~\ref{nans.thm.hh-digraph} yields that either
  $D / \theta(D) \cong k \cdot C_3^\circ + l \cdot \1^\circ$ for some $k, l \ge 0$ such that
  $k + l \ge 1$, or $D / \theta(D)$ is a finite homomorphism-homogeneous partially ordered set.
  Therefore, either $D$ is an inflation of
  $k \cdot C_3^\circ + l \cdot \1^\circ$ for some $k, l \ge 0$ such that $k + l \ge 1$, or
  $D$ is a finite homomorphism-homogeneous quasiorder (Corollary~\ref{nans.cor.quasiorder}).

  $(\Leftarrow)$
  Assume that $D$ belongs to one of the classes (1)--(2). Then $D / \theta(D)$
  is a homomorphism-homogeneous reflexive proper digraph according to Corollary~\ref{nans.cor.quasiorder}
  and Theorem~\ref{nans.thm.hh-digraph}.
  Lemma~\ref{nans.lem.al-retr} now yields that $D$ is a homomorphism-homogeneous improper digraph.
\end{PROOF}

The classification of bidirectionally disconnected digraphs with back-and-forth is slightly more
involved. Our intention is to prove that if $D$ is a homomorphism-homogeneous
bidirectionally disconnected digraph with back-and-forth, then $D / \theta(D)$ is
a homomorphism-homogeneous digraph with involution.

\begin{LEM}
  Let $D$ be a finite homomorphism-homogeneous reflexive bidirectionally disconnected improper digraph
  with back-and-forth.
  Then $S \sim T$ implies $S \rightleftarrows T$ for all $S$, $T \in V(D) / \theta(D)$.
\end{LEM}
\begin{PROOF}
  Take $S$, $T \in V(D) / \theta(D)$ so that $S \sim T$ and assume that $s \to t$ for some
  $s \in S$ and $t \in T$. We know that there exist distinct $U$, $W \in V(D) / \theta(D)$ such that
  $U \rightleftarrows W$, so choose $u_1, u_2 \in U$ and $w_1, w_2 \in W$ in such a way that
  $u_1 \to w_1$ and $w_2 \to u_2$. The mapping
  $$
    f : \begin{pmatrix}
      u_1 & w_1 \\
      s   & t
    \end{pmatrix}
  $$
  is a homomorphism between finitely induced subdigraphs of $D$, so it
  extends to an endomorphism $f^*$ of~$D$.
  It follows from Lemma~\ref{nans.lem.fS-S1} that $f^*(U) \subseteq S$ and $f^*(W) \subseteq T$, so
  $T \to S$ since $T \ni f^*(w_2) \to f^*(u_2) \in S$. Therefore, $S \rightleftarrows T$.
\end{PROOF}

\begin{LEM}\label{nans.lem.sTSt}
  Let $D$ be a finite homomorphism-homogeneous reflexive bidirectionally disconnected improper digraph,
  and let $S$ and $T$ be distinct classes of $\theta(D)$ such that $S \rightleftarrows T$.
  \begin{enumerate}
  \item
    There exists
    an $s \in S$ such that $s \rightleftarrows T$, and a $t \in T$ such that $S \rightleftarrows t$.
  \item
    $|S| \ge 2$ and $|T| \ge 2$.
  \item
    Suppose that $r \to t \to s$ for some $r, s \in S$ and $t \in T$. Then for every $u \in T$,
    either $r \to u \to s$ or $s \to u \to r$.
  \item
    $s \sim t$ for all $s \in S$ and all $t \in T$.
  \end{enumerate}
\end{LEM}
\begin{PROOF}
  Let $S$, $T \in V(D) / \theta(D)$ be distinct classes of $\theta(D)$ such that
  $S \rightleftarrows T$.

  (1)
  Take $s_1, s_2 \in S$ and $t_1, t_2 \in T$ so that $s_1 \to t_1$ and $t_2 \to s_2$.
  It suffices to show that there exists
  a $s \in S$ such that $s \rightleftarrows T$ \textbf{or} a $t \in T$ such that $S \rightleftarrows t$,
  since the mapping
  $$
   \begin{pmatrix}
     s_1 & t_1 \\
     t_2 & s_2
   \end{pmatrix}
  $$
  extends to an endomorphism od $D$ which then takes care of the other case.

  If $t_1 \to s_2$ then $S \rightleftarrows t_1$ and we are done. Assume now that $t_1 \not\to s_2$.
  Then the mapping
  $$
   f : \begin{pmatrix}
     s_1 & s_2 & t_1 \\
     s_1 & s_1 & t_1
   \end{pmatrix}
  $$
  is a homomorphism between finitely induced substructures of $D$ and, by the homogeneity requirement,
  extends to an endomorphism $f^*$ of $D$. From $t_2 \rightleftarrows t_1$ it follows that
  $f^*(t_2) \rightleftarrows t_1$, so $f^*(t_2) \in T$. Moreover, $t_2 \to s_2$ yields
  $f^*(t_2) \to s_1$. Therefore, $f^*(t_2) \to s_1 \to t_1$ and thus $s_1 \rightleftarrows T$.

  (2)
  Follows straightforwardly from (1) and the fact that $S$ and $T$ are disjoint classes of
  $\theta(D)$, so $s \not\rightleftarrows t$ for all $s \in S$ and all $t \in T$.

  (3)
  The statement trivially holds for $t$.
  Take any $u \in T \setminus \{t\}$ and let us show
  that the following two mappings cannot be homomorphisms between the corresponding induced substructures:
  $$
    f : \begin{pmatrix}
      r & s & u \\
      r & r & t
    \end{pmatrix}
    \text{\quad and\quad}
    g : \begin{pmatrix}
      s & r & u \\
      s & s & t
    \end{pmatrix}.
  $$
  Assume that $f$ is a homomorphism from $D[r, s, u]$ to $D[r, t]$. Then $f$ extends to an
  endomorphism $f^*$ of $D$. Let us take a look at $f^*(t)$. From $u \rightleftarrows t$ we infer
  $t = f^*(u) \rightleftarrows f^*(t)$, so $f^*(t) \in T$. On the other hand, $r \to t \to s$ implies
  $r = f^*(r) \to f^*(t) \to f^*(s) = r$, that is $r \rightleftarrows f^*(t)$, whence $f^*(t) \in S$.
  This contradicts the fact that $S$ and $T$ are disjoint.
  The proof for $g$ is analogous.

  Let us now show that $s \sim u$. Suppose, to the contrary, that $s \not\sim u$. If $u \not\to r$ then
  $f$ above is a homomorphism between finitely induced substructures of $D$, which is impossible.
  If, however, $u \to r$ then $g$ above is a homomorphism between finitely induced substructures of $D$,
  which is also impossible. Therefore, $s \sim u$.

  If $s \to u$ then $u \to r$ (since $s \to u$ and $u \not\to r$ implies that $f$ is a
  homomorphism between finitely induced substructures of $D$, which is impossible), and
  if $u \to s$ then $r \to u$ (since $u \to s$ and $r \not\to u$ implies that $g$ is a
  homomorphism between finitely induced substructures of $D$, which is impossible).

  (4)
  From (1) we know that there exist $q, r \in S$ and a $u \in T$ such that $q \to u \to r$.
  Take any $s \in S$ and any $t \in T$. Then (3) yields that $r \to t \to q$ or $q \to t \to r$.
  Clearly, if $s = r$ we are done, so we can assume that $s \ne r$.

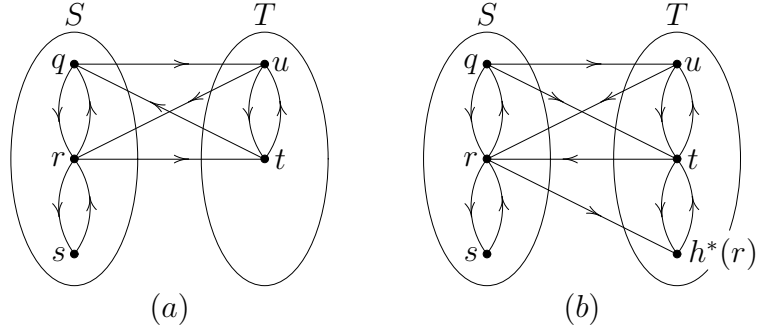
\begin{figure}
  \centering
  \unitlength 0.24pt
\begin{picture}(1200.0,500.0)
    \put(400.0,462.0){\makebox(0,0)[b]{$T$}}
    \put(100.0,462.0){\makebox(0,0)[b]{$S$}}
    \put(259.875,380.249){\arc{90.0}{0.58353}{1.10713}}
    \put(300.125,299.751){\arc{90.0}{4.24872}{4.77232}}
    \put(240.125,380.249){\arc{90.0}{1.51083}{2.03442}}
    \put(199.875,299.751){\arc{90.0}{5.17602}{5.69962}}
    \put(280.0,355.0){\arc{90.0}{4.18877}{4.71237}}
    \put(280.0,445.0){\arc{90.0}{1.57078}{2.09438}}
    \put(280.0,205.0){\arc{90.0}{4.18877}{4.71237}}
    \put(280.0,295.0){\arc{90.0}{1.57078}{2.09438}}
    \put(400.0,250.0){\ellipse{200.0}{400.0}}
    \put(100.0,250.0){\ellipse{200.0}{400.0}}
    \path(400.0,250.0)(100.0,400.0)
    \path(100.0,250.0)(400.0,250.0)
    \path(400.0,400.0)(100.0,250.0)
    \path(100.0,400.0)(400.0,400.0)
    \put(41.2235,154.452){\arc{70.2011}{5.4832}{6.15447}}
    \put(110.844,163.462){\arc{70.2011}{3.01287}{3.68414}}
    \put(158.777,195.546){\arc{70.2011}{2.34161}{3.01288}}
    \put(89.1561,186.537){\arc{70.2011}{6.15448}{6.82573}}
    \put(41.2235,304.452){\arc{70.2011}{5.4832}{6.15447}}
    \put(110.844,313.462){\arc{70.2011}{3.01287}{3.68414}}
    \put(158.777,345.546){\arc{70.2011}{2.34161}{3.01288}}
    \put(89.1561,336.537){\arc{70.2011}{6.15448}{6.82573}}
    \put(458.777,345.546){\arc{70.2011}{2.34161}{3.01288}}
    \put(389.156,336.537){\arc{70.2011}{6.15448}{6.82573}}
    \put(341.224,304.452){\arc{70.2011}{5.4832}{6.15447}}
    \put(410.844,313.462){\arc{70.2011}{3.01287}{3.68414}}
    \put(425.0,250.0){\makebox(0,0){$t$}}
    \put(425.0,400.0){\makebox(0,0){$u$}}
    \put(400.0,400.0){\circle*{12.0}}
    \put(400.0,250.0){\circle*{12.0}}
    \put(500.0,325.0){\arc{250.0}{2.49807}{3.78507}}
    \put(300.0,325.0){\arc{250.0}{5.63966}{6.92665}}
    \put(75.0,100.0){\makebox(0,0){$s$}}
    \put(75.0,250.0){\makebox(0,0){$r$}}
    \put(75.0,400.0){\makebox(0,0){$q$}}
    \put(0.0,325.0){\arc{250.0}{5.63966}{6.92665}}
    \put(0.0,175.0){\arc{250.0}{5.63966}{6.92665}}
    \put(200.0,175.0){\arc{250.0}{2.49807}{3.78507}}
    \put(200.0,325.0){\arc{250.0}{2.49807}{3.78507}}
    \put(100.0,100.0){\circle*{12.0}}
    \put(100.0,250.0){\circle*{12.0}}
    \put(100.0,400.0){\circle*{12.0}}
    \put(909.875,119.751){\arc{90.0}{4.65242}{5.17602}}
    \put(950.125,200.249){\arc{90.0}{2.03442}{2.55802}}
    \path(750.0,250.0)(1050.0,100.0)
    \put(1068.0,100.0){\makebox(0,0)[l]{$h^*(r)$}}
    \put(1050.0,100.0){\circle*{12.0}}
    \put(1150.0,175.0){\arc{250.0}{2.49807}{3.78507}}
    \put(950.0,175.0){\arc{250.0}{5.63966}{6.92665}}
    \put(1039.16,186.537){\arc{70.2011}{6.15448}{6.82573}}
    \put(1108.78,195.546){\arc{70.2011}{2.34161}{3.01288}}
    \put(1060.84,163.462){\arc{70.2011}{3.01287}{3.68414}}
    \put(991.224,154.452){\arc{70.2011}{5.4832}{6.15447}}
    \put(870.0,295.0){\arc{90.0}{1.04718}{1.57078}}
    \put(870.0,205.0){\arc{90.0}{4.71237}{5.23597}}
    \put(849.875,299.751){\arc{90.0}{4.65242}{5.17602}}
    \put(890.125,380.249){\arc{90.0}{2.03442}{2.55802}}
    \put(1050.0,462.0){\makebox(0,0)[b]{$T$}}
    \put(750.0,462.0){\makebox(0,0)[b]{$S$}}
    \put(909.875,380.249){\arc{90.0}{0.58353}{1.10713}}
    \put(950.125,299.751){\arc{90.0}{4.24872}{4.77232}}
    \put(930.0,355.0){\arc{90.0}{4.18877}{4.71237}}
    \put(930.0,445.0){\arc{90.0}{1.57078}{2.09438}}
    \put(1050.0,250.0){\ellipsearc{200.0}{400.0}{-5.21777}{0.597387}}
    \put(750.0,250.0){\ellipse{200.0}{400.0}}
    \path(1050.0,250.0)(750.0,400.0)
    \path(750.0,250.0)(1050.0,250.0)
    \path(1050.0,400.0)(750.0,250.0)
    \path(750.0,400.0)(1050.0,400.0)
    \put(691.224,154.452){\arc{70.2011}{5.4832}{6.15447}}
    \put(760.844,163.462){\arc{70.2011}{3.01287}{3.68414}}
    \put(808.777,195.546){\arc{70.2011}{2.34161}{3.01288}}
    \put(739.156,186.537){\arc{70.2011}{6.15448}{6.82573}}
    \put(691.224,304.452){\arc{70.2011}{5.4832}{6.15447}}
    \put(760.844,313.462){\arc{70.2011}{3.01287}{3.68414}}
    \put(808.777,345.546){\arc{70.2011}{2.34161}{3.01288}}
    \put(739.156,336.537){\arc{70.2011}{6.15448}{6.82573}}
    \put(1108.78,345.546){\arc{70.2011}{2.34161}{3.01288}}
    \put(1039.16,336.537){\arc{70.2011}{6.15448}{6.82573}}
    \put(991.224,304.452){\arc{70.2011}{5.4832}{6.15447}}
    \put(1060.84,313.462){\arc{70.2011}{3.01287}{3.68414}}
    \put(1075.0,250.0){\makebox(0,0){$t$}}
    \put(1075.0,400.0){\makebox(0,0){$u$}}
    \put(1050.0,400.0){\circle*{12.0}}
    \put(1050.0,250.0){\circle*{12.0}}
    \put(1150.0,325.0){\arc{250.0}{2.49807}{3.78507}}
    \put(950.0,325.0){\arc{250.0}{5.63966}{6.92665}}
    \put(725.0,100.0){\makebox(0,0){$s$}}
    \put(725.0,250.0){\makebox(0,0){$r$}}
    \put(725.0,400.0){\makebox(0,0){$q$}}
    \put(650.0,325.0){\arc{250.0}{5.63966}{6.92665}}
    \put(650.0,175.0){\arc{250.0}{5.63966}{6.92665}}
    \put(850.0,175.0){\arc{250.0}{2.49807}{3.78507}}
    \put(850.0,325.0){\arc{250.0}{2.49807}{3.78507}}
    \put(750.0,100.0){\circle*{12.0}}
    \put(750.0,250.0){\circle*{12.0}}
    \put(750.0,400.0){\circle*{12.0}}
    \put(250.0,38.0){\makebox(0,0)[t]{$(a)$}}
    \put(900.0,38.0){\makebox(0,0)[t]{$(b)$}}
\end{picture}
  \caption{The proof of Lemma \ref{nans.lem.sTSt} (4)}
  \label{nans.fig.1b1}
\end{figure}
  Assume, first, that $r \to t \to q$, Fig.~\ref{nans.fig.1b1}~$(a)$. Then, clearly, $t \ne u$.
  Since $u \to r \to t$, then from (3) we know that either $u \to s \to t$ or $t \to s \to u$.
  Either way, $s \sim t$.

  Assume, now, that $q \to t \to r$, Fig.~\ref{nans.fig.1b1}~$(b)$. The mapping
  $$
    h : \begin{pmatrix}
      s & t \\
      t & r
    \end{pmatrix}
  $$
  is a homomorphism between finitely induced substructures of $D$, so
  it extends to an endomorphism $h^*$ of $D$. Then $r \rightleftarrows s$ implies $h^*(r) \rightleftarrows t$,
  so $h^*(r) \in T$. Moreover, $t \to r$ implies $r \to h^*(r)$. Thus we get $t \to r \to h^*(r)$,
  so (3) ensures that $t \to s \to h^*(r)$ or $h^*(r) \to s \to t$. Either way, $s \sim t$.
\end{PROOF}

Let $S$ and $T$ be distinct classes of $\theta(D)$ such that $S \rightleftarrows T$. Define
a binary relation $\gamma_T(S) \subseteq S^2$ on $S$ as follows:
$$
    (x, y) \in \gamma_T(S) \text{\quad if and only if\quad}
    \lnot \exists t \in T \: (x \to t \to y \lor y \to t \to x).
$$

\begin{LEM}\label{nans.lem.RST}
  Let $D$ be a finite homomorphism-homogeneous reflexive improper bidirectionally disconnected digraph,
  and let $S$ and $T$ be distinct classes of $\theta(D)$ such that $S \rightleftarrows T$.
  \begin{enumerate}
  \item
    $\gamma_T(S)$ is an equivalence relation on $S$.
  \item
    $(x, y) \in \gamma_T(S)$ if and only if $\forall t \in T \; (x \rightarrow t \leftarrow y \lor
    x \leftarrow t \rightarrow y)$.
  \item
    $(x, y) \in \gamma_T(S)$ if and only if $\exists t \in T \; (x \rightarrow t \leftarrow y \lor
    x \leftarrow t \rightarrow y)$.
  \item
    $\gamma_T(S)$ has precisely two blocks.
  \item
    Let $S / \gamma_T(S) = \{ S_1, S_2 \}$ and $T / \gamma_S(T) = \{ T_1, T_2 \}$. Then
    $S_1 \TO T_1 \TO S_2 \TO T_2 \TO S_1$ or
    $T_1 \TO S_1 \TO T_2 \TO S_2 \TO T_1$.
  \end{enumerate}
\end{LEM}
\begin{PROOF}
  (1)
  The relation $\gamma_T(S)$ is obviously reflexive (because $S$ and $T$ are distinct classes of
  $\theta(D)$) and symmetric. Let us show that $\gamma_T(S)$ is transitive. Take any
  $(q, r)$, $(r, s) \in S^2$ and assume that $(q, s) \not\in \gamma_T(S)$.
  Then there exists a $t \in T$ such that $q \to t \to s$ or $s \to t \to q$. Without loss of
  generality we can assume that $q \to t \to s$. From Lemma~\ref{nans.lem.sTSt}~(4) we know that
  $x \sim y$ for all $x \in S$ and all $y \in T$, so $r \sim t$. Then $r \to t$ implies $(r, s) \not\in \gamma_T(S)$,
  while $t \to r$ implies $(q, r) \not\in \gamma_T(S)$. This shows that $\gamma_T(S)$ is an
  equivalence relation on~$S$.

  (2)
  Direction from right to left is obvious. In order to show the other direction, take any
  $(r, s) \in \gamma_T(S)$ and any $t \in T$. From Lemma~\ref{nans.lem.sTSt}~(4) we know that
  $x \sim y$ for all $x \in S$ and all $y \in T$, so $r \sim t \sim s$. Since
  $\lnot(r \to t \to s)$ and $\lnot(s \to t \to r)$, it must be the case that
  $r \rightarrow t \leftarrow s$ or $r \leftarrow t \rightarrow s$.

  (3)
  Direction from left to right follows straightforwardly from (2).
  In order to show the other direction, take any $(r, s) \not\in \gamma_T(S)$. Then there is a $t \in T$
  satisfying $r \to t \to s$ or $s \to t \to r$. Without loss of generality we can assume that
  $r \to t \to s$. Then Lemma~\ref{nans.lem.sTSt}~(3) ensures that
  $r \to u \to s$ or $s \to u \to r$ for every $u \in T$, whence, using Lemma~\ref{nans.lem.sTSt}~(4),
  it follows that $\lnot \exists u \in T \; (r \rightarrow u \leftarrow s \lor
  r \leftarrow u \rightarrow s)$.

  (4)
  We have shown in Lemma~\ref{nans.lem.sTSt}~(1)
  that there exist $r, s \in S$ and a $t \in T$ such that $r \to t \to s$. Then $r / \gamma_T(S) \ne
  s / \gamma_T(S)$ and thus $\gamma_T(S)$ has at least two blocks. Take any $q \in S$.
  From Lemma~\ref{nans.lem.sTSt}~(4) it follows that $q \sim t$.
  From (3) we now easily infer that $q \to t$ implies $q / \gamma_T(S) = r / \gamma_T(S)$, while
  $t \to q$ implies $q / \gamma_T(S) = s / \gamma_T(S)$. Therefore, $\gamma_T(S)$ has precisely two blocks.

  (5)
  Take any $s \in S_1$ and any $t \in T_1$. Then $s \sim t$ according to Lemma~\ref{nans.lem.sTSt}~(4).
  Without loss of generality we can assume that $s \to t$.
  Let us show that $S_1 \TO T_1$. Take any $r \in S$ and any $u \in T$. From (2) we conclude that
  $r \to t$ since $(r, s) \in \gamma_T(S)$. From $r \to t$ and $(t, u) \in \gamma_S(T)$ we infer
  $r \to u$. Therefore, $S_1 \TO T_1$. Now, $T_1$ and $T_2$ are distinct classes of
  $\gamma_S(T)$, so $S_1 \TO T_1$ implies $T_2 \to S_1$, and using the same argument as above we can show that
  $T_2 \TO S_1$. Analogously, $T_1 \TO S_2$ and $S_2 \TO T_2$.
\end{PROOF}

\begin{LEM}\label{nans.lem.STU-indep}
  Let $D$ be a finite homomorphism-homogeneous reflexive improper bidirectionally disconnected digraph,
  and let $S$, $T$ and $U$ be three distinct classes of $\theta(D)$ such that $S \rightleftarrows T$ and
  $T \rightleftarrows U$. Then $\gamma_S(T) = \gamma_U(T)$, that is,
  $\gamma_S(T)$ does not depend on~$S$.
\end{LEM}
\begin{PROOF}
  Assume, to the contrary, that $\gamma_S(T) \ne \gamma_U(T)$ and that there exists a pair
  $(t, t') \in \gamma_S(T)$ such that $(t, t') \notin \gamma_U(T)$. Then the definition of $\gamma$
  and Lemma~\ref{nans.lem.RST} provide us with an $s \in S$ and a $u \in U$ so that
  $t \rightarrow s \leftarrow t'$ and $t \to u \to t'$. The mapping
  $$
    f : \begin{pmatrix}
      s & t & u\\
      s & t & s
    \end{pmatrix}
  $$
  is a homomorphism between finitely induced substructures of $D$, so
  it extends to an endomorphism $f^*$ of $D$. From $t' \rightleftarrows t$ it follows that
  $f^*(t') \rightleftarrows t$, so $f^*(t') \in T$. On the other hand,
  $u \to t' \to s$ implies $f^*(t') \rightleftarrows s$, so $f^*(t') \in S$. Contradiction.
\end{PROOF}

\begin{THM}\label{nans.thm.case1b}
  Let $D$ be a finite reflexive bidirectionally disconnected improper digraph
  with back-and-forth.
  Then $D$ is homomorphism-homogeneous if and only if $D$ is an inflation of a
  homomorphism-homogeneous digraph with involution.
\end{THM}
\begin{PROOF}
  $(\Rightarrow)$
  Assume that $D$ is homomorphism-homogeneous.
  For a class $S$ of $\theta(D)$, define $\gamma(S)$ as follows:
  \begin{itemize}
  \item
    if $S \rightleftarrows T$ for some class $T$ of $\theta(D)$ distinct from $S$, let $\gamma(S) = \gamma_T(S)$;
  \item
    if $S \rightleftarrows T$ for no class $T$ of $\theta(D)$ distinct from $S$, let $\gamma(S) = S^2$.
  \end{itemize}
  With $\gamma(S)$ defined for every $S \in V(D) / \theta(D)$, we define $\gamma(D)$ by
  $$
    \gamma(D) = \UNION_{S \in V(D) / \theta(D)} \gamma(S).
  $$
  Then $D / \gamma(D)$ is well defined (Lemmas~\ref{nans.lem.STU-indep} and~\ref{nans.lem.RST})
  and it is a retract of $D$ (Lemma~\ref{nans.lem.al-retr}), so $D / \gamma(D)$ is homomorphism-homogeneous.
  Moreover, from Lemma~\ref{nans.lem.RST}~(4) and~(5), $D / \gamma(D)$ is a digraph with involution since
  every $\theta(D)$-class consists of at most two $\gamma(D)$-classes.

  $(\Leftarrow)$
  Let $D$ be an inflation of a homomorphism-homogeneous digraph with involution.
  According to Lemma~\ref{hhdl.lem.disconn-iff} it suffices to show that $D[S] \Rightarrow D[S']$ for all
  connected components $S$, $S'$ of $D$. Let $S$ and $S'$ be connected components of $D$ and let
  $f : U \to W$ be a homomorphism from $D[U]$ to $D[W]$ where $U \subseteq S$ and $W \subseteq S'$.
  Let $D = R\langle V_1, \ldots, V_n \rangle$ be an inflation of $R$ where $R$ is one of the digraphs
  listed in the statement of Theorem~\ref{nans.thm.hh-dwi},
  and let $\gamma(D)$ be the equivalence relation on $V(D)$ whose blocks are
  $V_1$, \ldots, $V_n$, so that $D / \gamma(D) \cong R$. Therefore, $R$ is a retract of $D$, so
  there exists a retraction-projection pair $r : V(D) \to V(R)$ and $j : V(R) \to V(D)$
  such that $r \circ j = \id$.

  Assume, first, that for every $\gamma(D)$-class $F$ there exists a $\gamma(D)$-class $F'$ such that
  $f(U \sec F) \subseteq F'$.
  Then $g : U / \gamma(D) \to W / \gamma(D)$ defined by $g(u / \gamma(D)) = f(u) / \gamma(D)$ is well-defined
  and it is a homomorphism from $R[U / \gamma(D)]$ to $R[W / \gamma(D)]$.
  We know that $R$ is homomorphism-homogeneous, so there exists a homomorphism
  $g^*$ from $R[S / \gamma(D)]$ to $R[S' / \gamma(D)]$ which extends $g$. But then
  the mapping $f^* : S \to S'$ defined by
  $$
    f^*(x) = \begin{cases}
      f(x), & x \in U\\
      j \circ g^* \circ r(x), & x \in S \setminus U
    \end{cases}
  $$
  is a homomorphism from $D[S]$ to $D[S']$ which extends $f$.

  Assume, now, that there exists a $\gamma(D)$-class $F$ such that $f(U \sec F)$ spreads over at least two
  $\gamma(D)$-classes. Then $f(U \sec F)$ spreads over exactly two $\gamma(D)$-classes $F_1$, $F_2$
  which belong to the same $\theta(D)$-class $H$, that is, $H = F_1 \union F_2$. Choose $x_1$, $x_2 \in U \sec F$
  in such a way that $f(x_1) \in F_1$ and $f(x_2) \in F_2$. Let us show that $f(U) \subseteq H$.
  Clearly, if $T$ is a $\theta(D)$-class that contains $F$ then $f(T \sec U) \subseteq H$.
  Let, now, $T$ be a $\theta(D)$-class such that $T \sec U \ne 0$ and $F \sec T = \0$.
  Take any $y \in T \sec U$. Since $D$ is an inflation of $R$
  it follows that $y \TO \{x_1, x_2\}$ or $\{x_1, x_2\} \TO y$, say $y \TO \{x_1, x_2\}$.
  Then $f(y) \TO \{f(x_1), f(x_2)\}$. Since $f(x_1)$ and $f(x_2)$ belong to distinct $\gamma(D)$-classes of the
  same $\theta(D)$-class, if $f(y) \notin H$ then $f(x_1) \to y$ and $f(x_2) \not\to y$, or $f(x_2) \to y$
  and $f(x_1) \not\to y$. Therefore, $f(y) \in H$. This shows that $f(U) \subseteq H$.
  Since $D[H] \cong K_{n_i}^\circ$, it is now easy to extend $f$ to a homomorphism $f^*$ from $D[S]$
  to $D[S']$: take any $h \in H$ and define $f^*$ by
  $$
    f^*(x) = \begin{cases}
      f(x), & x \in U\\
      h, & x \in S \setminus U.
    \end{cases}
  $$
  This concludes the proof.
\end{PROOF}

\begin{COR}\label{nans.cor.MAIN}
  Let $D$ be a finite reflexive binary relational system.
  If $D$ is bidirectionally disconnected, then $D$ is homomorphism-homogeneous if and only if
    \begin{enumerate}
    \item
      $D$ is a finite homomorphism-homogeneous quasiorder; or
    \item
      $D$ is an inflation of $k \cdot C_3^\circ + l \cdot \1^\circ$ for some $k, l \ge 0$ such that
      $k + l \ge 1$; or
    \item
      $D$ is an inflation of a homomorphism-homogeneous digraph with involution.
    \end{enumerate}

  If $D$ is bidirectionally connected then the problem of deciding whether $D$ is homomorphism-homogeneous
  is coNP-complete.
\end{COR}

\end{document}